\documentclass[11pt]{article}
\usepackage{amsmath}
\usepackage{amssymb}
\def\R{\mathbb{R}}
\def\N{\mathbb{N}}

\def\K{\mathop{\mathbf{K}}}
\newtheorem{theo}{\hspace*{\parindent}Theorem}
\newtheorem{lemma}{\hspace*{\parindent}Lemma}
\newtheorem{corol}{\hspace*{\parindent}Corollary}

\newtheorem{conjecture}{\hspace*{\parindent}Conjecture}
\newcounter{theremark}

\title{Log-convexity and log-concavity for series in gamma ratios and applications}
\author{S.I.\,Kalmykov and D.B.\,Karp}
\date{}
\begin{document}
\maketitle

\begin{center}
\parbox{12cm}{
\small\textbf{Abstract.} Polynomial sequence $\{P_m\}_{m\geq0}$ is $q$-logarithmically
concave if $P_{m}^2-P_{m+1}P_{m-1}$ is a polynomial with nonnegative coefficients
for any $m\geq{1}$.  We introduce an analogue of this notion for formal power series whose
coefficients are nonnegative continuous functions of  parameter.  Four types of such
power series are considered where parameter dependence is expressed by a ratio of gamma
functions.  We prove six theorems stating various forms of
$q$-logarithmic concavity and convexity of these series. The main motivating examples for these investigations
are hypergeometric functions. In the last section of the paper we present new
inequalities for the Kummer function, the ratio of the Gauss functions and the
generalized hypergeometric function obtained as direct applications of the
general theorems.

}
\end{center}

\bigskip

Keywords: \emph{Gamma function, log-concavity, log-convexity, q-log-concavity,
Wright log-concavity, Tur\'{a}n inequality, Kummer function, generalized hypergeometric function}

\bigskip

MSC2010: 26A51, 33C20, 33C15, 33C05

\bigskip

\paragraph{1. Introduction.} We adopt standard notation $\N$ for the set of positive integers,
$\N_{0}:=\N\cup\{0\}$, $\R$ will stand for reals and $\R_{+}$ for nonnegative reals.
The gamma function $\Gamma(x)$ was introduced by Leonard Euler who also
demonstrated that its second logarithmic derivative  is positive for positive values
of $x$.  In modern language  this means that $\Gamma(x)$ is
logarithmically convex (i.e. its logarithm is a convex function).
A sum of log-convex functions can be shown to be log-convex using
H\"{o}lder inequality or a theorem of Montel \cite[Theorem~1.4.5.2]{Mitr}.
Additivity implies then that the (finite or infinite) sum $f(\mu;x):=\sum f_k\Gamma(\mu+k)x^k$ is logarithmically
convex function of $\mu$ for fixed $x\geq{0}$ once the coefficients $f_k$ are assumed to be nonnegative.
It is not difficult to see that much more is true \cite[Theorem~2]{KS}: the formal
power series $f(\mu;x)f(\mu+\alpha+\beta;x)-f(\mu+\alpha;x)f(\mu+\beta;x)$
has nonnegative coefficients at all powers of $x$ if $\alpha,\beta\geq{0}$.
In \cite{KK} we considered a similar problem for the series $g(\mu;x):=\sum g_k\{\Gamma(\mu+k)\}^{-1}x^k$.
Here each term is log-concave function of $\mu$, so that lack of additivity of logarithmic
concavity does not allow to draw any immediate conclusion about the sum.
We have demonstrated, however, that the sequence $\{g(\mu;x)\}_{\mu\in\N}$
is log-concave for fixed $x>0$ if the sequence of coefficients $\{g_k\}_{k\in\N}$ is log-concave.
Moreover, in this case $g(\mu;x)g(\mu+\alpha+\beta;x)-g(\mu+\alpha;x)g(\mu+\beta;x)$ has
nonnegative coefficients at all powers of $x$ if $\alpha,\beta\in\N$.
The two sums above can be generalized naturally to series in product
ratios of gamma function having the form (\ref{eq:problem1}) below. Several known questions in financial mathematics \cite{BGR,CG},
multidimensional statistics \cite{SK}, probability \cite{NP} and special functions
\cite{BariczExpo,BarIsm,IsLaf} reduce to or depend on log-convexity or log-concavity
of special cases of such generalized series. Similar coefficient-wise positivity of product differences
is also important in combinatorics.  The following definition is
attributed to Richard Stanley \cite[p.795]{Sagan}.  A sequence of polynomials
$\{P_m(x)\}_{m\geq0}$ is said to be $q$-log-concave if
$$
P_{m}(x)^2-P_{m+1}(x)P_{m-1}(x)
$$
is a polynomial with nonnegative coefficients for any $m\geq{1}$.
It is strongly $q$-log-concave if
$$
P_{m}(x)P_{n}(x)-P_{m+1}(x)P_{n-1}(x)
$$
is a polynomial with nonnegative coefficients for all $m\geq{n}\geq{1}$.
The latter  notion was introduced by Sagan \cite{Sagan}. Many  sequences
of combinatorial polynomials especially those related to $q$-calculus possess
these properties (see \cite{CWY,Sagan} for details and references).
We will need extensions of these notions
to families of formal power series.  To be consistent with
the standard definitions of log-concavity and Wright log-concavity
\cite[Chapter~I.4]{MPF}, \cite[Section~1.1]{PPT} and to make our
formulations more compact,  we found it reasonable to change the
combinatorial terminology slightly.  Suppose
\begin{equation}\label{eq:fgeneral}
f(\mu;x)=\sum\limits_{k=0}^{\infty}f_k(\mu)x^k
\end{equation}
is a formal power series with nonnegative coefficients which depend continuously
on a nonnegative parameter $\mu$.

\noindent\textbf{Definition.} The family $\{f(\mu;x)\}_{\mu\geq{0}}$
is \textbf{Wright $q$-log-concave} if formal power series
\begin{equation}\label{eq:phigeneral}
\phi_{\mu}(\alpha,\beta;x):=f(\mu+\alpha;x)f(\mu+\beta;x)-f(\mu;x)f(\mu+\alpha+\beta;x)
\end{equation}
has nonnegative coefficients at all powers of $x$ for all $\mu,\alpha,\beta\geq{0}$.
If this property only holds
for $\alpha\in\N$ and all $\mu,\beta\geq{0}$ we will say that
$\{f(\mu;x)\}_{\mu\geq{0}}$ is \textbf{discrete Wright $q$-log-concave}.
Finally, $\{f(\mu;x)\}_{\mu\geq{0}}$ is \textbf{discrete $q$-log-concave} if
$\phi_{\mu}(\alpha,\beta;x)$ has nonnegative coefficients at all powers of $x$
for $\alpha\in\N$, $\beta\geq{\alpha-1}$ and all $\mu\geq{0}$.

If each function $f:\R_+\to\R_+$ is associated with the family of formal power series
$\{f(\mu;x)\}_{\mu\geq0}$ with constant term $f_0=f(\mu)$ and zero coefficients at all positive
powers of $x$, the above definitions become consistent with the following
standard terminology: $\mu\to{f(\mu)}$ is called Wright log-concave on $\R_+$
if $f(\mu+\alpha)f(\mu+\beta)\geq{f(\mu)f(\mu+\alpha+\beta)}$ for all $\mu,\alpha,\beta\geq{0}$
\cite[Chapter~I.4]{MPF}, \cite[Definition 1.13]{PPT}; it is discrete Wright
log-concave on $\R_+$ if the above inequality holds
for $\alpha\in\N$ and all $\mu,\beta\geq{0}$ and discrete
log-concave  if it holds for $\alpha\in\N$, $\beta\geq{\alpha-1}$
and $\mu\geq{0}$ \cite{KK}. For continuous functions Wright log-concavity is
equivalent to log-concavity (i.e. concavity of the logarithm). Discrete
Wright log-concavity implies discrete log-concavity but not vice versa
(see details in \cite{KK}). All above definitions also apply if
we substitute ''concave" by ''convex", ''non-negative" by
''non-positive" and reverse the sign of all inequalities.
In the theory of special functions discrete log-concavity and log-convexity are also
frequently referred to as ''Tur\'{a}n type inequalities'' following
the classical result of Paul Tur\'{a}n for Legendre polynomials
\cite{Turan}: $[P_n(x)]^2>P_{n-1}(x)P_{n+1}(x)$, $-1<x<1$.
Note, however, that the sequence $\{P_n(x)\}_{n\geq0}$ is not $q$-log-concave.
General Wright convex functions attracted a lot of attention recently
(see, for instance, \cite{IR,MP} and references therein)
following a fundamental result of Ng \cite{Ng}.

If $f:\N_0\to{\R_+}$ is a sequence, then discrete log-concavity
reduces to inequality $f_k^2\geq{f_{k-1}f_{k+1}}$, $k\in\N$.
We will additionally require that the sequence
$\{f_k\}_{k=0}^{\infty}$ is non-trivial and has no internal zeros, i.e.
$f_{N}=0$ implies either $f_{N+i}=0$ for all $i\in\N_{0}$ or
$f_{N-i}=0$ for $i=0,1,\ldots,N$.  Such
sequences are also known as $PF_2$ (P\'{o}lya frequency sub two)
or doubly positive \cite{Karlin}.  Clearly, if $f$ is (Wright)
log-concave then $1/f$ is (Wright) log-convex. Notwithstanding
the simplicity of this relation, several important properties
of log-concavity and log-convexity differ. As we already mentioned above,
log-convexity is preserved under addition while log-concavity is
not. Further, log-convexity is a stronger property than convexity
whereas log-concavity is weaker than concavity.  Further properties of
log-convex and log-concave functions can be found, for instance, in
\cite[3E, 16D, 18B]{MOA}, \cite[Chapter~2]{NP}  and  \cite[Chapter~13]{PPT}.

The questions considered in \cite{KK,KS} and in this paper are particular cases of the following
general problem: under what conditions on a nonnegative sequence $\{f_k\}$ and
the numbers $a_i$, $b_j$ the series
\begin{equation}\label{eq:problem1}
f(\mu;x)=\sum\limits_{k=0}^{\infty}f_k\frac{\prod_{i=1}^{n}\Gamma(a_i+\mu+\varepsilon_ik)}
{\prod_{j=1}^{m}\Gamma(b_j+\mu+\varepsilon_{n+j}k)}x^k
\end{equation}
is (discrete, Wright) $q$-log-concave or $q$-log-convex?
Here $\varepsilon_r$ can take values $1$ or $0$. In particular,
if the ratio $f_{k+1}/f_{k}$ is a
rational function of $k$ the series in (\ref{eq:problem1}) is
hypergeometric (possibly times some gamma functions) and $\mu$ represents
parameter shift \cite[Chapter~2]{AAR}.

In \cite{KS} the authors treated the cases of (\ref{eq:problem1}) with $n=1$,
$m=0$, $\varepsilon_1=1$; $n=m=1$, $\varepsilon_1=1$,
$\varepsilon_2=0$, $a_1=b_1$; and $n=m=1$, $\varepsilon_1=0$,
$\varepsilon_2=1$, $a_1=b_1$.  In \cite{KK} we handled $n=0$, $m=1$, $\varepsilon_1=1$.
In this paper we treat the following
cases:

(a) $n=m=2$, $\varepsilon_1=1$, $\varepsilon_2=0$, $\varepsilon_3=0$, $\varepsilon_4=1$,  $a_1=b_1$, $a_2=b_2$;

(b) $n=m=1$, $\varepsilon_1=1$, $\varepsilon_2=1$;

(c) $n=2$, $m=1$, $\varepsilon_1=1$, $\varepsilon_2=0$, $\varepsilon_3=1$, $a_2=b_1$;

(d) $n=1$, $m=2$, $\varepsilon_1=1$, $\varepsilon_2=0$, $\varepsilon_3=1$, $a_1=b_1$.

For small values of $n$ and $m$ considered in (a)-(d) it is easy to determine
the conditions for each term in (\ref{eq:problem1}) to be
log-convex. By additivity we can then assert the log-convexity of $\mu\to{f(\mu;x)}$
for fixed $x\geq0$, but not $q$-log-convexity
of any type, i.e. non-positivity of the Taylor coefficients of $\phi_{\mu}(\alpha,\beta;x)$
defined by (\ref{eq:phigeneral}). If the terms in (\ref{eq:problem1}) are
log-concave even the verification of log-concavity of $\mu\to{f(\mu;x)}$ for fixed $x\geq0$
becomes non-trivial. In this paper we verify Wright $q$-log-convexity
and discrete $q$-log-concavity for the family of power series defined
in (\ref{eq:problem1}) under restrictions (a)-(d).
Our results will imply then that either $x\to\phi_{\mu}(\alpha,\beta;x)$ or
$x\to-\phi_{\mu}(\alpha,\beta;x)$ is absolutely monotonic.
According to Hardy, Littlewood and P\'{o}lya theorem
\cite[Proposition~2.3.3]{NP} absolute monotonicity of  $x\to\phi_{\mu}(\alpha,\beta;x)$ implies that
this function is multiplicatively convex:
$$
\phi_{\mu}(\alpha,\beta;x^{\lambda}y^{1-\lambda})\leq{\phi_{\mu}(\alpha,\beta;x)^{\lambda}\phi_{\mu}(\alpha,\beta;y)^{1-\lambda}}
$$
for $\lambda\in[0,1]$. Curiously, this inequality leads to interesting results
even when applied to the simplest function $1+x^2$. We have
$$
1+x^{2\lambda}y^{2(1-\lambda)}\leq (1+x^2)^{\lambda}(1+y^2)^{1-\lambda}~~\text{and}~~1+x^{2(1-\lambda)}y^{\lambda}\leq (1+x^2)^{1-\lambda}(1+y^2)^{\lambda}.
$$
Multiplying these inequalities and simplifying we obtain
$$
(x^{\lambda}y^{1-\lambda})^2+(x^{1-\lambda}y^{\lambda})^2\leq x^2+y^2,
$$
which is equivalent to inequality
$$
M_2(G_{\lambda}(x,y),G_{1-\lambda}(x,y))\leq M_2(x,y),~~x,y\geq0, ~~0\leq\lambda\leq1,
$$
where $M_2(a,b)=\sqrt{(a^2+b^2)/2}$ is quadratic mean,
$G_{\lambda}(a,b)=a^{\lambda}b^{1-\lambda}$ is weighted geometric mean.

The paper is organized as follows: in section~2 we collect
several lemmas repeatedly used in the proofs; section~3 comprises
six theorems constituting the main content of the paper; in
section~4 we present several applications and relate them to
some known results.

\paragraph{2. Preliminaries.}
We will need several lemmas which we prove in this section.
\begin{lemma}\label{lm:product}
Suppose $\{f(\mu;x)\}_{\mu\geq0}$ and $\{g(\mu;x)\}_{\mu\geq0}$
are \emph{(}discrete, Wright\emph{)} $q$-log-concave. Then $\{f(\mu;x)g(\mu;x)\}_{\mu\geq0}$
is  \emph{(}discrete, Wright\emph{)}  $q$-log-concave.
\end{lemma}
\textbf{Remark.} Lemma~\ref{lm:product} holds, in particular, if $\mu\to{g(\mu)}$ is
a log-concave function independent of $x$.
\textbf{Proof.} We have
\begin{multline*}
f(\mu+\alpha;x)g(\mu+\alpha;x)f(\mu+\beta;x)g(\mu+\beta;x)-f(\mu;x)g(\mu;x)f(\mu+\alpha+\beta;x)g(\mu+\alpha+\beta;x)=
\\
g(\mu+\alpha;x)g(\mu+\beta;x)(f(\mu+\alpha;x)f(\mu+\beta;x)-f(\mu;x)f(\mu+\alpha+\beta;x))
\\
+f(\mu;x)f(\mu+\alpha+\beta;x)(g(\mu+\alpha;x)g(\mu+\beta;x)-g(\mu;x)g(\mu+\alpha+\beta;x)).
\end{multline*}
This formula implies the claim of the lemma.~$\square$

\begin{lemma}\label{lm:disc-Wright}
Let $f$ be a nonnegative-valued function defined on $\R_+$. Suppose
$$
\phi_{\mu}(\alpha,\beta):=f(\mu+\alpha)f(\mu+\beta)-f(\mu)f(\mu+\beta+\alpha)\geq{0}~\text{for}~\alpha=1~\text{and all}~\mu,\beta\geq{0}.
$$
Then $\phi_{\mu}(\alpha,\beta)\geq{0}$ for all $\alpha\in\N$ and $\mu,\beta\geq{0}$,
i.e. $\mu\to{f(\mu)}$ is discrete Wright log-concave on $\R_+$.
\end{lemma}
\textbf{Proof.}  According to assumptions of the lemma written
for the pairs $\{\mu,\beta\}$, $\{\mu+1,\beta\}$, $\{\mu,\beta+1\}$ and $\alpha=1$
we have
\begin{gather}
f(\mu+1)f(\mu+\beta)\geq{f(\mu)f(\mu+\beta+1),}\label{eq:Wright1}
\\[5pt]
f(\mu+2)f(\mu+\beta+1)\geq{f(\mu+1)f(\mu+\beta+2),}\label{eq:Wright2}
\\[5pt]
f(\mu+1)f(\mu+\beta+1)\geq{f(\mu)f(\mu+\beta+2).}\label{eq:Wright3}
\end{gather}
A product of (\ref{eq:Wright1}) and (\ref{eq:Wright2}) reads
$$
f(\mu+1)f(\mu+\beta+1)(f(\mu+\beta)f(\mu+2)-f(\mu)f(\mu+\beta+2))\geq0.
$$
This implies either $f(\mu+\beta)f(\mu+2)\geq{f(\mu)f(\mu+\beta+2)}$
which is our claim for $\alpha=2$ or $f(\mu+1)f(\mu+\beta+1)=0$ which implies
$f(\mu)f(\mu+\beta+2)=0$ according to (\ref{eq:Wright3}), so that again
$f(\mu+\beta)f(\mu+2)\geq{f(\mu)f(\mu+\beta+2)}$.
Hence, we have demonstrated that $\phi_{\mu}(2,\beta)\geq{0}$.
In a similar fashion $\phi_{\mu}(\alpha,\beta)\geq{0}$ holds for
all $\alpha\in\N$ and $\mu,\beta\geq{0}$.~$\square$

In the above Lemma the function $f$ may or may not be defined by the series
(\ref{eq:fgeneral}) - we have not made any use of the special series structure
 in the proof. In the next lemma we deal with Wright $q$-log-concavity
 and the series definition becomes important.

\begin{lemma}\label{lm:disc-qWright}
Let the series $f$ be defined by  \emph{(\ref{eq:fgeneral})} and suppose
$$
\phi_{\mu}(1,\beta;x):=f(\mu+1;x)f(\mu+\beta;x)-f(\mu;x)f(\mu+\beta+1;x)
$$
has nonnegative coefficients at powers of $x$ for and all $\mu,\beta\geq{0}$.
Then $\phi_{\mu}(\alpha,\beta;x)$ has nonnegative coefficients
at powers of $x$ for all $\alpha\in\N$, $\beta\geq\alpha-1$ and $\mu\geq{0}$.
\end{lemma}
\textbf{Proof.}  Define
$$
\psi_{\nu,x}(a,b):=f(\nu;x)^2-f(\nu-a;x)f(\nu+b;x).
$$
By assumptions of the lemma the difference
\begin{multline*}
\psi_{\nu,x}(a,b)-\psi_{\nu,x}(a-1,b-1)=f(\nu-a+1;x)f(\nu+b-1;x)
-f(\nu-a;x)f(\nu+b;x)\\
\overset{\normalfont\mu:=\nu-a}{=}f(\mu+1;x)f(\mu+a+b-1;x)-f(\mu;x)f(\mu+a+b;x)
\end{multline*}
has nonnegative power series coefficients when $\nu\geq{a}$,
$a+b-1\geq{0}$.  Further, for a positive integer $k$
\begin{multline*}
f(\mu+k;x)f(\mu+\beta;x)-f(\mu;x)f(\mu+\beta+k;x)
=f(\nu-a+k;x)f(\nu+b-k;x)-f(\nu-a;x)f(\nu+b;x)
\\[3pt]
=\psi_{\nu,x}(a,b)-\psi_{\nu,x}(a-k,b-k)=
(\psi_{\nu,x}(a,b)-\psi_{\nu,x}(a-1,b-1))+
\\[3pt]
(\psi_{\nu,x}(a-1,b-1)-\psi_{\nu,x}(a-2,b-2))+\cdots
+(\psi_{\nu,x}(a-k+1,b-k+1)-\psi_{\nu,x}(a-k,b-k)),
\end{multline*}
where $\mu=\nu-a$, $\mu+\beta=\nu+b-k$. We must require
$\nu\geq{a}$, $a+b-1\geq{0}$ for the first parentheses to have nonnegative
power series coefficients, $\nu\geq{a-1}$, $a+b-3\geq{0}$ for the second parentheses
to have nonnegative power series coefficients, and so on up to
$\nu\geq{a-k+1}$, $a+b-2k+1\geq{0}$.  These inequalities reduce to
$\mu\geq{0}$ and $\beta\geq{k-1}$.~$\square$

The next lemma is implied by a much stronger result of Alzer
\cite[Theorem~10]{Alzer}.
\begin{lemma}\label{lm:Gammaratio}
Suppose $0\leq\min(\alpha_1,\alpha_2)\leq\min(\beta_1,\beta_2)$
and $\alpha_1+\alpha_2\leq\beta_1+\beta_2$. Then the function
\[
x\to\frac{\Gamma(x+\alpha_1)\Gamma(x+\alpha_2)}{\Gamma(x+\beta_1)\Gamma(x+\beta_2)}
\]
is strictly monotone decreasing on $(0,\infty)$, except when the sets $\{\alpha_1,\alpha_2\}$
and $\{\beta_1,\beta_2\}$ are equal.
\end{lemma}
Next, we formulate an elementary inequality we will repeatedly use
below.
\begin{lemma}\label{lm:uvrs}
Suppose $u,v,r,s>0$, $u=\max(u,v,r,s)$ and $uv>rs$. Then
$u+v>r+s$.
\end{lemma}

Lemma~\ref{lm:uvrs} is a particular case of a much more general
result on logarithmic majorization - see \cite[2.A.b]{MOA}.
See also \cite[Lemma~1]{KK} for a direct proof.

In the next lemma we say that a sequence has no more than one
change of sign if it has the pattern $(--\cdots--00\cdots00++\cdots++)$,
where zeros and minus signs may be missing.
\begin{lemma}\label{lm:sum}
Suppose $\{f_k\}_{k=0}^{n}$ has no internal zeros and
$f_k^2\geq{f_{k-1}f_{k+1}}$, $k=1,2,\ldots,n-1$. If the real
sequence $A_0,A_1,\ldots,A_{[n/2]}$ satisfying $A_{[n/2]}>0$ and
$\sum\limits_{0\leq{k}\leq{n/2}}\!\!\!\!A_k\geq{0}$ has no more
than one change of sign, then
\begin{equation}\label{eq:keysum}
\sum\limits_{0\leq{k}\leq{n/2}}f_{k}f_{n-k}A_k\geq{0}.
\end{equation}
Equality is only attained if $f_k=f_0\alpha^k$, $\alpha>0$, and
$\sum\limits_{0\leq{k}\leq{n/2}}\!\!\!\!A_k=0$.
\end{lemma}
A proof of this lemma is found in \cite[Lemma~2]{KK}.

The generalized hypergeometric
function is defined by the series
\begin{equation}\label{eq:pFqdefined}
{_{p}F_q}\left(\left.\!\!\begin{array}{c} A\\
B\end{array}\right|z\!\right)={_{p}F_q}\left(A;B;z\right):=\sum\limits_{n=0}^{\infty}\frac{(a_1)_n(a_2)_n\cdots(a_{p})_n}{(b_1)_n(b_2)_n\cdots(b_q)_nn!}z^n,
\end{equation}
where we write $A=(a_1,a_2,\ldots,a_p)$, $B=(b_1,b_2,\ldots,b_q)$
for brevity and $(a)_0=1$, $(a)_n=a(a+1)\cdots(a+n-1)$,
$n\geq{1}$,  denotes the rising factorial. The series
(\ref{eq:pFqdefined}) converges in the entire complex plane if
$p\leq{q}$ and in the unit disk if $p=q+1$.  In the latter case
its sum can be extended analytically to the whole complex plane
cut along the ray $[1,\infty)$ \cite[Chapter~2]{AAR}.  The
series (\ref{eq:pFqdefined}) is a particular case of
(\ref{eq:problem1}) because $(a)_k=\Gamma(a+k)/\Gamma(a)$.

The next identity for the Kummer function ${_1F_1}$ is believed to be new
and may be of independent interest.
\begin{lemma}\label{lm:Kummer-id}
The Kummer function ${_1F_1}$ satisfies the following identity\emph{:}
\begin{multline}\label{eq:Kummer-id}
{_1F_1}(a+\mu;c+\mu;x){_1F_1}(a+1;c+1;x)-{_1F_1}(a+\mu+1;c+\mu+1;x){_1F_1}(a;c;x)\\=
\frac{(c-a)x}{c(c+1)(c+\mu)(c+\mu+1)}\left[(c+\mu)(c+\mu+1){_1F_1}(a+1;c+2;x){_1F_1}(a+\mu+1;c+\mu+1;x)\right.
\\\left.-c(c+1){_1F_1}(a+1;c+1;x){_1F_1}(a+\mu+1;c+\mu+2;x)\right].
\end{multline}
\end{lemma}
\textbf{Proof.}  Apply the easily verifiable contiguous relations
\[
{_1F_1}(a;c;x)={_1F_1}(a;c+1;x)+\frac{ax}{c(c+1)}{_1F_1}(a+1;c+2;x),
\]
\[
{_1F_1}(a+\mu;c+\mu;x)={_1F_1}(a+\mu+1;c+\mu+1;x)-\frac{(c-a)x}{(c+\mu)(c+\mu+1)}{_1F_1}(a+\mu+1;c+\mu+2;x)
\]
and
\[
{_1F_1}(a+1;c+1;x)={_1F_1}(a;c+1;x)+\frac{x}{c+1}{_1F_1}(a+1;c+2;x)
\]
to the corresponding functions on the left hand side of
(\ref{eq:Kummer-id}). Expanding and collecting similar terms we
can then rewrite the left-hand side of (\ref{eq:Kummer-id}) as
\begin{multline*}
\frac{(c-a)x}{c(c+1)(c+\mu)(c+\mu+1)}\left[(c+\mu)(c+\mu+1){_1F_1}(a+1;c+2;x){_1F_1}(a+\mu+1;c+\mu+1;x)\right.
\\\left.-c((c+1){_1F_1}(a;c+1;x)+x{_1F_1}(a+1;c+2;x)){_1F_1}(a+\mu+1;c+\mu+2;x)\right].
\end{multline*}
Finally, applying here the contiguous relation
\[
(c+1){_1F_1}(a;c+1;x)+x{_1F_1}(a+1;c+2;x)=(c+1){_1F_1}(a+1;c+1;x),
\]
yields the right hand-side of
(\ref{eq:Kummer-id}).~$\square$

The next lemma has been proved using some ideas borrowed from
\cite{CHJ}.
\begin{lemma}\label{lm:ab-sum}
The inequality
\begin{equation}\label{eq:sum}
\sum\limits_{k=0}^{m}\frac{(a)_k(a+\mu)_{m-k}}{(b)_k(b+\mu)_{m-k}}\binom{m}{k}(m-2k+\mu)\geq{0},
\end{equation}
holds for each integer $m\geq{1}$ and all $\mu\geq{0}$ if
$b\geq{a}\geq{0}$ or $a\geq{b}\geq{1}$.
\end{lemma}
\textbf{Proof.} Denote
$$
u_k=\frac{(a)_k(a+\mu)_{m-k}}{(b)_k(b+\mu)_{m-k}}.
$$
If $a=b$ or $a=0$ the claim is obvious. Suppose first that
$b>a>0$. Then $x\to(a+x)/(b+x)$ increasing so that for $k<m-k$
$$
u_k>u_{m-k},~~\text{since}~~\frac{(a+\mu+k)\cdots(a+\mu+m-k-1)}{(b+\mu+k)\cdots(b+\mu+m-k-1)}>\frac{(a+k)\cdots(a+m-k-1)}{(b+k)\cdots(b+m-k-1)}.
$$
It follows that
$$
\binom{m}{k}u_k(m-2k+\mu)+\binom{m}{m-k}u_{m-k}(2k-m+\mu)=\binom{m}{k}[(m-2k)(u_k-u_{m-k})+\mu(u_k+u_{m-k})]>0
$$
for each $k\leq{m-k}$ which proves the lemma for all
$b\geq{a}\geq{0}$.  If $a>b\geq{1}$ things become more
complicated. In this case we will apply Abel's lemma (summation by
parts) in the form \cite{CHJ}
$$
\sum\limits_{k=0}^{m}(\alpha_{k+1}-\alpha_{k})\beta_{k}
=\sum\limits_{k=0}^{m}\alpha_{k+1}(\beta_{k}-\beta_{k+1})+\alpha_{m+1}\beta_{m+1}-\alpha_{0}\beta_{0}.
$$
Gosper's algorithm \cite{Gosper}, \cite[Chapter~5]{PWZ} produces the following
antidifference which is easy to verify directly:
$$
u_k(m-2k+\mu)=\alpha_{k+1}-\alpha_{k},~\text{where}~\alpha_k=\frac{(b-1)(b-1+\mu)(a)_k(a+\mu)_{m+1-k}}{(a-b+1)(b-1)_{k}(b-1+\mu)_{m+1-k}},~k=0,1,\ldots,m+1.
$$
Next, setting $\beta_k=\binom{m}{k}$ we immediately obtain
$$
\beta_{k}-\beta_{k+1}=\binom{m}{k}-\binom{m}{k+1}=\frac{2k+1-m}{m+1}\binom{m+1}{k+1}.
$$
Hence, an application of Abel's lemma yields (we use the fact that $\beta_{-1}=\beta_{m+1}=0$):
\begin{multline*}
\sum\limits_{k=0}^{m}\binom{m}{k}u_k(m-2k+\mu)=\sum\limits_{k=0}^{m}\alpha_{k+1}(\beta_{k}-\beta_{k+1})-\alpha_{0}\beta_{0}
\\
=\frac{(b-1)(b-1+\mu)}{(a-b+1)}\left\{\sum\limits_{k=0}^{m}\frac{(a)_{k+1}(a+\mu)_{m-k}}{(b-1)_{k+1}(b-1+\mu)_{m-k}}
\frac{2k+1-m}{m+1}\binom{m+1}{k+1}-\frac{(a+\mu)_{m+1}}{(b-1+\mu)_{m+1}}\right\}
\\
=\frac{(b-1)(b-1+\mu)}{(a-b+1)}\left\{\sum\limits_{k=-1}^{m}\frac{(a)_{k+1}(a+\mu)_{m-k}}{(b-1)_{k+1}(b-1+\mu)_{m-k}}
\frac{2k+1-m}{m+1}\binom{m+1}{k+1}\right\}
\\
=\frac{(b-1)(b-1+\mu)}{(a-b+1)n}\left\{\sum\limits_{j=0}^{n}\frac{(a)_{j}(a+\mu)_{n-j}}{(b-1)_{j}(b-1+\mu)_{n-j}}
\binom{n}{j}(2j-n)\right\},
\end{multline*}
where $j=k+1$, $n=m+1$. If $a>b\geq{1}$ and $0\leq{j}<n-j$ we have
$$
\frac{(a)_{j}(a+\mu)_{n-j}}{(b-1)_{j}(b-1+\mu)_{n-j}}<\frac{(a)_{n-j}(a+\mu)_{j}}{(b-1)_{n-j}(b-1+\mu)_{j}}
$$$$
~\Leftrightarrow~
\frac{(a+\mu+j)\cdots(a+\mu+n-j-1)}{(b-1+\mu+j)\cdots(b+\mu+n-j-2)}<\frac{(a+j)\cdots(a+n-j-1)}{(b-1+j)\cdots(b+n-j-2)}
$$
since $x\to(a+x)/(b+x)$ is decreasing.  Hence,
\begin{multline*}
\sum\limits_{j=0}^{n}\frac{(a)_{j}(a+\mu)_{n-j}}{(b-1)_{j}(b-1+\mu)_{n-j}}\binom{n}{j}(2j-n)
\\
=\sum\limits_{0\leq{j}<n/2}
\left(\frac{(a)_{n-j}(a+\mu)_{j}}{(b-1)_{n-j}(b-1+\mu)_{j}}-\frac{(a)_{j}(a+\mu)_{n-j}}{(b-1)_{j}(b-1+\mu)_{n-j}}\right)\binom{n}{j}(n-2j)>0.~\square
\end{multline*}

\noindent\textbf{Remark.} With more effort one can show that (\ref{eq:sum})
remains valid if $a\geq{b}\geq{1/2}$, but we will not use this
fact in the present paper.

\paragraph{3. Main results.} In this section we prove six general
theorems for series in  ratios of
rising factorials and gamma functions.  The power series
expansions in this section  are understood as formal, so that no
questions of convergence are discussed. In applications the radii
of convergence will usually be apparent. The results of this
section are exemplified by concrete special functions in the
subsequent section. The first two theorems deal with the class of series  defined by
\begin{equation}\label{eq:f-def}
f_{a,c}(\mu;x):=\sum\limits_{n=0}^{\infty}f_n\frac{(a+\mu)_n}{(c+\mu)_n}\frac{x^n}{n!}.
\end{equation}
Since
\begin{multline*}
f_{a,c}(\mu+\alpha;x)f_{a,c}(\mu+\beta;x)-f_{a,c}(\mu;x)f_{a,c}(\mu+\alpha+\beta;x)
\\
=f_{a+\mu,c+\mu}(\alpha;x)f_{a+\mu,c+\mu}(\beta;x)-f_{a+\mu,c+\mu}(0;x)f_{a+\mu,c+\mu}(\alpha+\beta;x),
\end{multline*}
there is no loss of generality in considering the product difference (\ref{eq:phigeneral})
in the form
\begin{equation}\label{eq:phi-def}
\phi_{a,c}(\mu,\nu;x):=f_{a,c}(\mu;x)f_{a,c}(\nu;x)-f_{a,c}(0;x)f_{a,c}(\mu+\nu;x)=\sum\limits_{m=1}^{\infty}\phi_mx^m.
\end{equation}
Logarithmic concavity or convexity of $\mu\mapsto
f_{a,c}(\mu;x)$ depends on the interrelation between $a$ and $c$.
\begin{theo}\label{th:f-concavity}
Suppose $c\geq{a}>0$ and $\{f_n\}_{n=0}^{\infty}$ is a nonnegative
log-concave sequence without internal zeros. Then the
function $\mu\mapsto{f_{a,c}(\mu;x)}$ is discrete
Wright log-concave on $[0,\infty)$ for each fixed $x>0$. Moreover, the family
$\{f_{a,c}(\mu;x)\}_{\mu\geq0}$ is discrete $q$-log-concave,
i.e. the function $x\to\phi_{a,c}(\mu,\nu;x)$ has nonnegative
power series coefficients for all $\nu\in\N$ and $\mu\geq\nu-1$ so that
$x\to\phi_{a,c}(\mu,\nu;x)$ is absolutely monotonic and
multiplicatively convex on $(0,\infty)$.
\end{theo}
\textbf{Remark.} It easy to see from the proof of the theorem
that $\phi_m>0$ for all $m\geq{1}$ if $f_n>0$ for all $n\geq{0}$, $c>a$
and $\mu>0$.

\textbf{Proof.} If $c=a$ the claim is obvious. Suppose $c>a>0$.
 According to Lemmas~\ref{lm:disc-Wright} and \ref{lm:disc-qWright}
it is sufficient to consider the case $\nu=1$. For a fixed integer $m\geq{1}$ we have by the Cauchy
product and Gauss summation:
$$
\phi_m\!=\!\sum\limits_{k=0}^{m}\!f_{k}f_{m-k}\!\underbrace{\left(\frac{(a+1)_k(a+\mu)_{m-k}}{(c+1)_k(c+\mu)_{m-k}k!(m-k)!}-
\frac{(a)_k(a+\mu+1)_{m-k}}{(c)_k(c+\mu+1)_{m-k}k!(m-k)!}\right)}_{N_k}\!\!=\!\!
\sum\limits_{k=0}^{[m/2]}f_{k}f_{m-k}M_k,
$$
where $M_k=N_{k}+N_{m-k}$ for $k<m/2$ and $M_k=N_{k}$ for $k=m/2$.
We aim to apply Lemma~\ref{lm:sum} to prove that $\phi_m\geq0$.
First, we need to show that
\begin{equation}\label{eq:sumNk}
\sum\limits_{k=0}^{[m/2]}M_k=\sum\limits_{k=0}^{m}N_k\geq0.
\end{equation}
Since, clearly,
\[
\sum\limits_{m=1}^{\infty}x^m\sum\limits_{k=0}^{m}N_k=
{_1F_1}(a+\mu;c+\mu;x){_1F_1}(a+1;c+1;x)-{_1F_1}(a+\mu+1;c+\mu+1;x){_1F_1}(a;c;x),
\]
we are in the position to apply formula (\ref{eq:Kummer-id}) from
Lemma~\ref{lm:Kummer-id} which, after equating power series coefficients
on both sides, yields:
\begin{multline*}
\sum\limits_{k=0}^{m+1}N_{k}=\frac{(c-a)}{c(c+1)(c+\mu)(c+\mu+1)}\times
\\
\sum\limits_{k=0}^{m}\Biggl(\underbrace{\frac{(c+\mu)(c+\mu+1)(a+1)_k(a+\mu+1)_{m-k}}{(c+2)_k(c+\mu+1)_{m-k}k!(m-k)!}}_{u_k}
-\underbrace{\frac{c(c+1)(a+1)_k(a+\mu+1)_{m-k}}{(c+1)_k(c+\mu+2)_{m-k}k!(m-k)!}}_{r_k}\Biggr)
\\=\sideset{}{_{}^{\,\prime}}{\sum}\limits_{k=0}^{[m/2]}(u_{k}+u_{k-m}-r_{k}-r_{k-m}).
\end{multline*}
Here the prime at the summation sign means that the term with
$k=m/2$ (which only happens for even $m$) has multiplier $1/2$.
This last term is positive since ($l=m/2$)
\[
u_{l}>r_{l}~\Leftrightarrow~\frac{(c+\mu)(c+\mu+1)}{(c+2)_l(c+\mu+1)_{l}}
>\frac{c(c+1)}{(c+1)_l(c+\mu+2)_{l}}~\Leftrightarrow~
\frac{(c+\mu)(c+\mu+l+1)}{c(c+l+1)}>1.
\]
We claim that all other terms in the rightmost sum above are also
positive by Lemma~\ref{lm:uvrs} with $u=u_{k}$, $v=u_{m-k}$,
$r=r_k$, $s=r_{m-k}$. To verify the assumptions of
Lemma~\ref{lm:uvrs} we need to show that $u_{k}>u_{k-m}$,
$u_{k}>r_{k}$, $u_{k}>r_{m-k}$ and $u_{k}u_{k-m}>r_{k}r_{k-m}$. We
have
$$
u_{k}>u_{k-m}~\Leftrightarrow~\left.\frac{\Gamma(a+1+x)\Gamma(c+\mu+1+x)}{\Gamma(c+2+x)\Gamma(a+\mu+1+x)}\right|_{x=k}
>\left.\frac{\Gamma(a+1+x)\Gamma(c+\mu+1+x)}{\Gamma(c+2+x)\Gamma(a+\mu+1+x)}\right|_{x=m-k},
$$
since the gamma quotient is decreasing by
Lemma~\ref{lm:Gammaratio} and $k<m-k$ by assumption.  Next,
\[
u_{k}>r_{k}~\Leftrightarrow~(c+\mu)(c+\mu+1+m-k)>c(c+1+k),
\]
which is true by $\mu>0$ and $k<m-k$. The inequality
$u_{k}>r_{m-k}$ reduces to
$$
\left.\frac{(c+\mu)(c+\mu+k+1)\Gamma(a+1+x)\Gamma(c+\mu+1+x)}{c(c+k+1)\Gamma(a+\mu+1+x)\Gamma(c+1+x)}\right|_{x=k}
>\left.\frac{\Gamma(a+1+x)\Gamma(c+\mu+1+x)}{\Gamma(a+\mu+1+x)\Gamma(c+1+x)}\right|_{x=m-k},
$$
which is true because the gamma quotient is decreasing by
Lemma~\ref{lm:Gammaratio} while
$(c+\mu)(c+\mu+k+1)/[c(c+k+1)]\geq{1}$.
Finally,
\[
u_{k}u_{k-m}>r_{k}r_{k-m}~\Leftrightarrow~(c+\mu)^2(c+\mu+k+1)(c+\mu+m-k+1)>c^2(c+k+1)(c+m-k+1),
\]
which proves inequality (\ref{eq:sumNk}).

Next, we need to demonstrate that the sequence
$\{M_k\}_{k=0}^{[m/2]}$ changes sign not more than once. To this
end introduce the following notation
$$
\widetilde{M}_k=k!(m-k)!M_k=\left\{\begin{array}{ll}
\underbrace{\frac{(a+1)_k(a+\mu)_{m-k}}{(c+1)_k(c+\mu)_{m-k}}}_{=\widetilde{u}_{k}}
+\underbrace{\frac{(a+1)_{m-k}(a+\mu)_{k}}{(c+1)_{m-k}(c+\mu)_{k}}}_{=\widetilde{u}_{m-k}}
\\[25pt]
-\underbrace{\frac{(a)_k(a+\mu+1)_{m-k}}{(c)_k(c+\mu+1)_{m-k}}}_{=\widetilde{r}_{k}}
-\underbrace{\frac{(a)_{m-k}(a+\mu+1)_{k}}{(c)_{m-k}(c+\mu+1)_{k}}}_{=\widetilde{r}_{m-k}},&
k<m/2
\\[25pt]
\frac{\displaystyle(a+1)_{m/2}(a+\mu)_{m/2}}{\displaystyle(c+1)_{m/2}(c+\mu)_{m/2}}-\frac{\displaystyle(a)_{m/2}(a+\mu+1)_{m/2}}
{\displaystyle(c)_{m/2}(c+\mu+1)_{m/2}},& k=m/2.
\end{array}\right.
$$
Suppose that $\widetilde{M}_k<0$ for some $0<k<m/2$ then we will show that
$\widetilde{M}_{k-1}<0$. Indeed, $\widetilde{M}_{k-1}$ can be
written in the following form
$$
\widetilde{M}_{k-1}=\underbrace{\frac{(c+k)(a+\mu+m-k)}{(a+k)(c+\mu+m-k)}}_{=I_1}\widetilde{u}_k
+
\underbrace{\frac{(a+1+m-k)(c+\mu+k-1)}{(c+1+m-k)(a+\mu+k-1)}}_{=I_2}\widetilde{u}_{m-k}-
$$
$$
-\underbrace{\frac{(c+k-1)(a+\mu+1+m-k)}{(a+k-1)(c+\mu+1+m-k)}}_{=I_3}\widetilde{r}_k
-
\underbrace{\frac{(a+m-k)(c+\mu+k)}{(c+m-k)(a+\mu+k)}}_{=I_4}\widetilde{r}_{m-k}=
$$
$$
=I_1(\widetilde{u}_{k}+\widetilde{u}_{m-k}-\widetilde{r}_{k}-\widetilde{r}_{m-k})
+ (I_1-I_3)(\widetilde{r}_k-\widetilde{u}_{m-k})+
(I_2-I_3)(\widetilde{u}_{m-k}-\widetilde{r}_{m-k})+
(I_1+I_2-I_3-I_4)\widetilde{r}_{m-k}.
$$
The first term is negative since $\widetilde{M}_{k}<0$. We will
show that all other terms are also negative.  The function
$x\mapsto\frac{\displaystyle\beta+x}{\displaystyle\alpha+x}$,
$\beta>\alpha$, is strictly decreasing on $(0,\infty)$ which leads to the
following inequalities
$$
I_1<I_3 \Leftrightarrow
\frac{(c+k)(a+\mu+m-k)}{(a+k)(c+\mu+m-k)}<\frac{(c+k-1)(a+\mu+1+m-k)}{(a+k-1)(c+\mu+1+m-k)},
$$
$$
I_2<I_3 \Leftrightarrow
\frac{(a+1+m-k)(c+\mu+k-1)}{(c+1+m-k)(a+\mu+k-1)}<\frac{(c+k-1)(a+\mu+1+m-k)}{(a+k-1)(c+\mu+1+m-k)},
$$
$$
I_4<I_2 \Leftrightarrow
\frac{(a+m-k)(c+\mu+k)}{(c+m-k)(a+\mu+k)}<\frac{(a+1+m-k)(c+\mu+k-1)}{(c+1+m-k)(a+\mu+k-1)},
$$
valid for $0<k<m/2$ and $\mu>0$. Hence,
$I_3=\max(I_1,I_2,I_3,I_4)$.  Further, $I_3I_4>I_1I_2$ is
equivalent to
$$
H_1(\mu):=\frac{(a+\mu+k-1)(c+\mu+k)(c+\mu+m-k)(a+\mu+1+m-k)}{(c+\mu+k-1)(a+\mu+k)(a+\mu+m-k)(a+\mu+1+m-k)}>H_1(0).
$$
We will show that $H_1(\mu)$ is increasing on $(0,\infty)$.
Indeed, by straightforward calculation
$$
\frac{d}{d\mu}\log(H_1(\mu))=(c-a)(H_2(z_2+1)-H_2(z_2)-(H_2(z_1+1)-H_2(z_1))),
$$
where $0\leq{z_1}=k-1<z_2=m-k$ and
$$
H_2(z)=\frac{1}{(a+\mu+z)(c+\mu+z)}
$$
is convex on $[0,\infty)$ implying
$H_2(z_2+1)-H_2(z_2)>H_2(z_1+1)-H_2(z_1)$. Thus we have proved
that $I_3I_4>I_1I_2$ so that by Lemma~\ref{lm:uvrs} we get
$I_1+I_2-I_3-I_4<0$.  To demonstrate that $\widetilde{M}_{k-1}<0$
it remains to show that $\widetilde{u}_{m-k}>\widetilde{r}_{m-k}$
and $\widetilde{r}_k>\widetilde{u}_{m-k}$.  We have
$$
\widetilde{u}_{m-k}> \widetilde{r}_{m-k} \Leftrightarrow
\frac{(a+1)_{m-k}(a+\mu)_{k}}{(c+1)_{m-k}(c+\mu)_{k}}
>\frac{(a)_{m-k}(a+\mu+1)_{k}}{(c)_{m-k}(c+\mu+1)_{k}}\Leftrightarrow
\frac{(a+m-k)c}{a(c+m-k)}>\frac{(a+\mu+k)(c+\mu)}{(a+\mu)(c+\mu+k)}.
$$
Since $\mu\mapsto
\frac{\displaystyle(a+\mu+k)(c+\mu)}{\displaystyle(a+\mu)(c+\mu+k)}$
is strictly decreasing on $[0,\infty)$
$$
\frac{(a+\mu+k)(c+\mu)}{(a+\mu)(c+\mu+k)}<\frac{(a+k)c}{a(c+k)}<\frac{(a+m-k)c}{a(c+m-k)},
$$
where the rightmost inequality clearly holds for $0<k<m/2$.
Finally, in order to show that
$\widetilde{r}_k>\widetilde{u}_{m-k}$ it suffices to prove that
$\widetilde{u}_{k}>\widetilde{r}_{m-k}$. Indeed, $\widetilde{u}_{m-k}\geq\widetilde{r}_k$
and the preceding inequality imply that $\widetilde{M}_{k}>0$ contradicting our hypothesis.
The validity of $\widetilde{u}_{k}>\widetilde{r}_{m-k}$ follows from
$$
\widetilde{u}_{k}>\widetilde{r}_{m-k}~\Leftrightarrow~
\frac{(a+1)_k(a+\mu)_{m-k}}{(c+1)_k(c+\mu)_{m-k}}>
\frac{(a)_{m-k}(a+\mu+1)_{k}}{(c)_{m-k}(c+\mu+1)_{k}}
~\Leftrightarrow~
$$
$$
~\Leftrightarrow~
\frac{c\Gamma(a+1+k)\Gamma(a+\mu+m-k)}{a\Gamma(a+m-k)\Gamma(a+\mu+1+k)}>
\frac{(c+\mu)\Gamma(c+1+k)\Gamma(c+\mu+m-k)}{(a+\mu)\Gamma(c+m-k)\Gamma(c+\mu+1+k)}.
$$
It is easy to see that conditions of Lemma~\ref{lm:Gammaratio} are
satisfied for the gamma ratio for all $0\leq{k}<m/2$ and $\mu>0$,
while clearly $c/a>(c+\mu)/(a+\mu)$.~~$\square$

\begin{corol}\label{cr:phi-complmon-conc}
Suppose $c>a>0$ and the series in $(\ref{eq:f-def})$
converges for all $x\geq{0}$. Then for all $\nu\in\N$ and
$\mu\geq\nu-1$ the function $y\to\phi_{a,c}(\mu,\nu;1/y)$
is completely monotonic and log-convex on $[0,\infty)$
and there exists a nonnegative measure $\tau$ supported on $[0,\infty)$
such that
$$
\phi_{a,c}(\mu,\nu;x)=\int\limits_{[0,\infty)}e^{-t/x}d\tau(t).
$$
\end{corol}
\textbf{Proof}.  According to \cite[Theorem~3]{MS} a convergent
series of completely monotonic functions with nonnegative coefficients is
again completely monotonic. This implies that
$y\to\phi_{a,c}(\mu,\nu;1/y)$ is completely monotonic, so that  the
above integral representation follows by Bernstein's theorem
\cite[Theorem~1.4]{SSV}.  Log-convexity follows from complete
monotonicity according to \cite[Exersice 2.1(6)]{NP}.~~$\square$

\begin{corol}\label{cr:phi-below-conc}
Under hypotheses and notation of Theorem~\ref{th:f-concavity} and
for all $\nu\in\N$, $\mu\geq{\nu-1}$ and $x\geq{0}$
$$
f_{a,c}(\mu;x)f_{a,c}(\nu;x)-f_{a,c}(0;x)f_{a,c}(\mu+\nu;x)\geq
\frac{f_0f_1x\mu\nu(c-a)(2c+\mu+\nu)}{c(c+\mu)(c+\nu)(c+\mu+\nu)}
$$
with equality only at $x=0$ if $c-a,\mu,\nu\ne0$.
\end{corol}
\textbf{Proof.} Indeed, the claimed inequality is just
$\phi_{a,c}(\mu,\nu;x)\geq\phi_1x$ which is true by
Theorem~\ref{th:f-concavity}.~~$\square$

There is virtually no doubt that the discrete $q$-log-concavity
demonstrated in Theorem~\ref{th:f-concavity} results from our method
of proof so that the adjective ''discrete'' is redundant. In other
words, we propose the following conjecture.

\begin{conjecture} The family $\{f(\mu;x)\}_{\mu\geq{0}}$ is
Wright $q$-log-concave for all $c\geq{a}>0$.
\end{conjecture}

Next theorem handles the case  $a\geq{c}>0$.
As it turns out frequently the log-convexity case is
simpler.
\begin{theo}\label{th:f-convexity}
Suppose $a\geq{c}>0$, $\{f_n\}_{n=0}^{\infty}$ is any nonnegative
sequence and the functions $f_{a,c}(\mu;x)$ and $\phi_{a,c}(\mu,\nu;x)$
are defined by $(\ref{eq:f-def})$ and $(\ref{eq:phi-def})$, respectively.
Then the function $\mu\mapsto{f_{a,c}(\mu;x)}$ is strictly log-convex
on $[0,\infty)$ for each fixed $x>0$. Moreover, the family
$\{f_{a,c}(\mu;x)\}_{\mu\geq0}$ is Wright $q$-log-convex, i.e. the function
$x\to\phi_{a,c}(\mu,\nu;x)$ has non-positive  power series
coefficients so that $x\to-\phi_{a,c}(\mu,\nu;x)$ is
absolutely monotonic and multiplicatively convex on $(0,\infty)$.
\end{theo}

\textbf{Proof.} If $a=c$ the claim is obvious. Suppose $a>c>0$.
Combining the Cauchy product with the Gauss summation as in the proof
of Theorem~\ref{th:f-concavity} the problem reduces to the inequality
\begin{equation}\label{eq:HmMk}
-\phi_m=\sum\limits_{k=0}^{[m/2]}\frac{f_kf_{m-k}}{k!(m-k)!}M_k>0,
\end{equation}
where
$$
M_k=\underbrace{\frac{(a)_k(a+\mu+\nu)_{m-k}}{(c)_k(c+\mu+\nu)_{m-k}}}_{=v}
+\underbrace{\frac{(a)_{m-k}(a+\mu+\nu)_{k}}{(c)_{m-k}(c+\mu+\nu)_{k}}}_{=u}
-\underbrace{\frac{(a+\mu)_k(a+\nu)_{m-k}}{(c+\mu)_k(c+\nu)_{m-k}}}_{=r}
-\underbrace{\frac{(a+\mu)_{m-k}(a+\nu)_{k}}{(c+\mu)_{m-k}(c+\nu)_{k}}}_{=s}
$$
for $k<m/2$ and
$$
M_k=\frac{(a)_k(a+\mu+\nu)_{m-k}}{(c)_k(c+\mu+\nu)_{m-k}}
-\frac{(a+\mu)_k(a+\nu)_{m-k}}{(c+\mu)_k(c+\nu)_{m-k}}~~\text{for}~k=m/2
$$
(this term is missing for odd values of $m$).
We will show that $M_k>0$ for each $k=1,2,\ldots,m/2$.
We will need the following fact \cite[Lemma~2 and Remark~7]{KS}: the
function
\[
x\mapsto\frac{(x+\alpha_1)(x+\alpha_2)}{(x+\beta_1)(x+\beta_2)},~~\alpha_1,\alpha_2,\beta_1,\beta_2\geq{0},
\]
is increasing on $[0,\infty)$ iff
\[
\frac{\beta_1\beta_2}{\alpha_1\alpha_2}\geq\frac{\beta_1+\beta_2}{\alpha_1+\alpha_2}\geq{1}.
\]
It follows that this function is bounded by $1$ for all
$x\geq{0}$ since its value at infinity is $1$.
If any of the inequalities above is strict the
function is strictly increasing.  This implies that $M_k>0$ if
$k=m/2$ for
\[
\frac{(a)_k(a+\mu+\nu)_{k}}{(a+\mu)_k(a+\nu)_{k}}
>\frac{(c)_k(c+\mu+\nu)_{k}}{(c+\mu)_k(c+\nu)_{k}}.
\]
Indeed, both sides of this inequality represent a product of the
terms
\[
\frac{(x+i)(x+\mu+\nu+i)}{(x+\mu+i)(x+\nu+i)}
\]
Since $(\mu+i)(\nu+i)>i(\mu+\nu+i)$ for each nonnegative integer
$i$, this function is increasing and its value at $x=a$ is greater
then its value at $x=c<a$.

Further, we will show that  $M_k>0$ for $0\le{k}<m/2$  using
Lemma~\ref{lm:uvrs}.  We have
$$
u>v~\Leftrightarrow~\frac{\Gamma(c+k)\Gamma(a+\mu+\nu+k)}{\Gamma(a+k)\Gamma(c+\mu+\nu+k)}>
\frac{\Gamma(c+m-k)\Gamma(a+\mu+\nu+m-k)}{\Gamma(a+m-k)\Gamma(c+\mu+\nu+m-k)}
$$
by Lemma~\ref{lm:Gammaratio} for $a>c>0$ and $k<m-k$;
$$
u>s~\Leftrightarrow~\frac{(c+\nu)_k(a+\mu+\nu)_k}{(a+\nu)_k(c+\mu+\nu)_k}>
\frac{(c)_{m-k}(a+\mu)_{m-k}}{(a)_{m-k}(c+\mu)_{m-k}}
$$
because
$$
\frac{(c+\nu+i)(a+\mu+\nu+i)}{(a+\nu+i)(c+\mu+\nu+i)}>
\frac{(c+i)(a+\mu+i)}{(a+i)(c+\mu+i)}~\text{for}~i=0,1,\ldots,k-1,
$$
and $0<(c+i)(a+\mu+i)/[(a+i)(c+\mu+i)]<1$ for $i=k,\ldots,m-k-1$
by the fact above;  Next, $u>r$ by exactly the same argument with
$\mu$ and $\nu$ interchanged; finally,
$$
uv>rs~\Leftrightarrow~\frac{(a)_k(a+\nu+\mu)_k}{(a+\nu)_k(a+\mu)_k}\times\frac{(a)_{m-k}(a+\nu+\mu)_{m-k}}{(a+\nu)_{m-k}(a+\mu)_{m-k}}>
\frac{(c)_k(c+\nu+\mu)_k}{(c+\nu)_k(c+\mu)_k}\times\frac{(c)_{m-k}(c+\nu+\mu)_{m-k}}{(c+\nu)_{m-k}(c+\mu)_{m-k}},
$$
where the first multiplier on the left is bigger than the first
multiplier on the right and the second multiplier on the left is
bigger than the second multiplier on the right due to the
increase of $x\to(x+i)(x+\nu+\mu+i)/[(a+\nu+i)(a+\mu+i)]$ for each
$i\geq{0}$.  Hence, by Lemma~\ref{lm:uvrs} $M_k=u+v-r-s>0$.
~~$\square$

The following two corollaries are similar to
Corollaries~\ref{cr:phi-complmon-conc} and~\ref{cr:phi-below-conc}.
\begin{corol}\label{cr:phi-below-conc1}
Under hypotheses and notation of Theorem~\ref{th:f-convexity} and for
all $\mu,\nu,x\geq{0}$
$$
f_{a,c}(0;x)f_{a,c}(\mu+\nu;x)-f_{a,c}(\mu;x)f_{a,c}(\nu;x)
\!\geq\frac{f_0f_1x\mu\nu(a-c)(2c+\mu+\nu)}{c(c+\mu)(c+\nu)(c+\mu+\nu)}
$$
with equality only at $x=0$ if $a-c,\mu,\nu\ne0$.
\end{corol}
\begin{corol}\label{cr:phi-complmon-conv}
Suppose $a>c>0$ and the series in $(\ref{eq:f-def})$
converges for all $x\geq{0}$. Then for all $\nu,\mu>0$  the function
$y\to-\phi_{a,c}(\mu,\nu;1/y)$
is completely monotonic and log-convex on $[0,\infty)$
and there exists a nonnegative measure $\tau$ supported on $[0,\infty)$
such that
$$
\phi_{a,c}(\mu,\nu;x)=-\int\limits_{[0,\infty)}e^{-t/x}d\tau(t).
$$
\end{corol}

\bigskip

The next two theorems deal with  the class of series  defined by
\begin{equation}\label{eq:g-def}
\mu\to
g_{a,c}(\mu;x)=\sum\limits_{n=0}^{\infty}g_n\frac{\Gamma(a+n+\mu)}{\Gamma(c+n+\mu)}\frac{x^n}{n!}
\end{equation}
and their product differences
\begin{equation}\label{eq:psi-def}
\psi_{a,c}(\mu,\nu;x)=g_{a,c}(\mu;x)g_{a,c}(\nu;x)-g_{a,c}(0;x)g_{a,c}(\mu+\nu;x)
=\sum\limits_{m=0}^{\infty}\psi_mx^m.
\end{equation}
If we set $g_n=f_n$ we get
$$
g_{a,c}(\mu;x)=\frac{\Gamma(a+\mu)}{\Gamma(c+\mu)}f_{a,c}(\mu;x)
$$
where $f_{a,c}(\mu;x)$ is defined by (\ref{eq:f-def}). It is then tempting to
derive the properties of $g_{a,c}(\mu;x)$ from Theorems~\ref{th:f-concavity}
and \ref{th:f-convexity} using Lemma~\ref{lm:product}.
However, when $a>c$ the gamma ratio in front of $f_{a,c}(\mu;x)$ is log-concave while $\mu\to{f_{a,c}(\mu;x)}$
is log-convex, so that Lemma~\ref{lm:product} cannot be applied.
Similar situation holds for $c\geq{a}$.

\begin{theo}\label{th:g-convexity}
Suppose $c\geq{a}>0$, $\{g_n\}_{n=0}^{\infty}$ is any nonnegative
sequence.  Then the function $\mu\to g_{a,c}(\mu;x)$ is Wright log-convex
on $[0,\infty)$ for each fixed $x\geq{0}$. Moreover, the family
$\{g_{a,c}(\mu;x)\}_{\mu\geq0}$ is Wright $q$-log-convex, i.e.
the function $x\to\psi_{a,c}(\mu,\nu;x)$ has non-positive power
series coefficients for all $\mu,\nu\geq{0}$ so that $x\to-\psi_{a,c}(\mu,\nu;x)$ is absolutely
monotonic  and multiplicatively convex on $(0,\infty)$.
\end{theo}
\noindent\textbf{Proof.}  Cauchy product and Gauss summation yield
\begin{multline*}
-\psi_m\!=\!\sum\limits_{k=0}^{m}\frac{g_{k}g_{m-k}}{k!(m-k)!}
\!\left\{\frac{\Gamma(a+k)\Gamma(a+\nu+\mu+m-k)}{\Gamma(c+k)\Gamma(c+\nu+\mu+m-k)}
-\frac{\Gamma(a+\nu+k)\Gamma(a+\mu+m-k)}{\Gamma(c+\nu+k)\Gamma(c+\mu+m-k)}\!\right\}\!
\\
=\!\!\sum\limits_{k=0}^{[m/2]}\frac{g_{k}g_{m-k}}{k!(m-k)!}M_k,
\end{multline*}
where
\begin{multline*}
M_k=\underbrace{\frac{\Gamma(a+k)\Gamma(a+\nu+\mu+m-k)}{\Gamma(c+k)\Gamma(c+\nu+\mu+m-k)}}_{=u}
+\underbrace{\frac{\Gamma(a+m-k)\Gamma(a+\nu+\mu+k)}{\Gamma(c+m-k)\Gamma(c+\nu+\mu+k)}}_{=v}
\\
-\underbrace{\frac{\Gamma(a+\nu+k)\Gamma(a+\mu+m-k)}{\Gamma(c+\nu+k)\Gamma(c+\mu+m-k)}}_{=r}
-\underbrace{\frac{\Gamma(a+\nu+m-k)\Gamma(a+\mu+k)}{\Gamma(c+\nu+m-k)\Gamma(c+\mu+k)}}_{=s}
\end{multline*}
for $k<m/2$, and
\[
M_k=\frac{\Gamma(a+k)\Gamma(a+\nu+\mu+m-k)}{\Gamma(c+k)\Gamma(c+\nu+\mu+m-k)}
-\frac{\Gamma(a+\nu+k)\Gamma(a+\mu+m-k)}{\Gamma(c+\nu+k)\Gamma(c+\mu+m-k)}~~\text{for}~~k=m/2
\]
(this term is missing for odd values of $m$). We aim to
demonstrate that $M_k>0$ using Lemma~\ref{lm:uvrs}. We have
$$
u>v~\Leftrightarrow~\frac{\Gamma(a+k)\Gamma(c+\nu+\mu+k)}{\Gamma(c+k)\Gamma(a+\nu+\mu+k)}
>\frac{\Gamma(a+m-k)(c+\nu+\mu+m-k)}{\Gamma(c+m-k)\Gamma(a+\nu+\mu+m-k)}.
$$
In view of $k<m-k$, the last inequality holds by
Lemma~\ref{lm:Gammaratio} with $\alpha_1=c+\nu+\mu$, $\alpha_2=a$,
$\beta_1=\max(c,a+\nu+\mu)$ è $\beta_2=\min(c,a+\nu+\mu)$,
$x_1=k$, $x_2=m-k$.  Next,
$$
u>r~\Leftrightarrow~\frac{\Gamma(a+k)\Gamma(a+\nu+\mu+m-k)}{\Gamma(a+\nu+k)\Gamma(a+\mu+m-k)}
>
\frac{\Gamma(c+k)(c+\nu+\mu+m-k)}{\Gamma(c+\nu+k)(c+\mu+m-k)}
$$
Setting $\alpha_1=\nu+\mu+m-k$, $\alpha_2=k$,
$\beta_1=\max(\nu+k,\mu+m-k)$, $\beta_2=\min(\nu+k,\mu+m-k)$,
$x_1=a$ and $x_2=c$ we get the last inequality by
Lemma~\ref{lm:Gammaratio} again. In a similar fashion one can
demonstrate that $u>s$.  Finally, $uv>rs$ by multiplication of the
following two inequalities
\[
\frac{\Gamma(a+k)\Gamma(a+\nu+\mu+k)}{\Gamma(a+\nu+k)\Gamma(a+\mu+k)}
>\frac{\Gamma(c+k)\Gamma(c+\nu+\mu+k)}{\Gamma(c+\nu+k)\Gamma(c+\mu+k)}
\]
and
\[
\frac{\Gamma(a+m-k)\Gamma(a+\nu+\mu+m-k)}{\Gamma(a+\mu+m-k)\Gamma(a+\nu+m-k)}
>
\frac{\Gamma(c+m-k)\Gamma(c+\nu+\mu+m-k)}{\Gamma(c+\mu+m-k)\Gamma(c+\nu+m-k)},
\]
each of them holds by Lemma~\ref{lm:Gammaratio}.~~$\square$

Again we have two corollaries similar to
Corollaries~\ref{cr:phi-complmon-conc} and~\ref{cr:phi-below-conc}.
\begin{corol}\label{cr:psi-below-conv}
Under hypotheses and notation of Theorem~\ref{th:g-convexity} and
for all $\mu,\nu,x\geq{0}$
$$
g_{a,c}(0;x)g_{a,c}(\mu+\nu;x)-g_{a,c}(\mu;x)g_{a,c}(\nu;x)
\!\geq\!
g_0^2\left\{\frac{\Gamma(a)\Gamma(a+\mu+\nu)}{\Gamma(c)\Gamma(c+\mu+\nu)}
-\frac{\Gamma(a+\nu)\Gamma(a+\mu)}{\Gamma(c+\nu)\Gamma(c+\mu)}\right\}
$$
with equality only at $x=0$ if $c-a,\mu,\nu\ne0$.
\end{corol}
\begin{corol}\label{cr:psi-complmon-conv}
Suppose $a>c>0$ and the series in $(\ref{eq:g-def})$ converges
for all $x\geq{0}$. Then  for all $\nu,\mu>0$ the function
$y\to-\psi_{a,c}(\mu,\nu;1/y)$
is completely monotonic and log-convex on $[0,\infty)$
and there exists a nonnegative measure $\tau$ supported on $[0,\infty)$
such that
$$
\psi_{a,c}(\mu,\nu;x)=-\int\limits_{[0,\infty)}e^{-t/x}d\tau(t).
$$
\end{corol}
Next, we can combine Theorem~\ref{th:f-concavity} and Theorem~\ref{th:g-convexity}
to get
\begin{corol}\label{cr:f-twosided}
Under hypotheses and notation of Theorem~\ref{th:f-concavity}
$$
\frac{(c+\mu)_{\nu}(a)_{\nu}}{(a+\mu)_{\nu}(c)_{\nu}}\leq\frac{f_{a,c}(0;x)f_{a,c}(\mu+\nu;x)}{f_{a,c}(\nu;x)f_{a,c}(\mu;x)}\leq1
$$
for all $\nu\in\N$, $\mu\geq{0}$ and $x\geq{0}$.
\end{corol}

\textbf{Proof.} The estimate from above is a restatement of
$\phi_{a,c}(\mu,\nu;x)\geq0$ valid by Theorem~\ref{th:f-concavity}.
To demonstrate the estimate from below set in Theorem~\ref{th:g-convexity}
$$
g_{a,c}(\mu;x)=\frac{\Gamma(a+\mu)}{\Gamma(c+\mu)}f_{a,c}(\mu;x).
$$
In view of $(a)_k=\Gamma(a+k)/\Gamma(a)$ the required inequality is
a restatement of $\psi_{a,c}(\mu,\nu;x)\leq0$.~$\square$

Further, combining Corollary~\ref{cr:phi-below-conc} with Corollary~\ref{cr:f-twosided}
we obtain the following two-sided bounds for the Tur\'{a}nian:
\begin{equation}\label{eq:f-Turanian}
\frac{2xf_0f_1\nu^2(c-a)}{c(c+\nu)(c+2\nu)}\leq f_{a,c}(\nu;x)^2-f_{a,c}(0;x)f_{a,c}(2\nu;x)\leq\frac{(a+\nu)_{\nu}(c)_{\nu}-(c+\nu)_{\nu}(a)_{\nu}}{(c)_{\nu}(a+\nu)_{\nu}}f_{a,c}(\nu;x)^2.
\end{equation}
This holds for $\nu\in\N$, $c\geq{a}>0$,  $x\geq{0}$ and a
log-concave sequence $\{f_n\}_{n\geq0}$ without internal zeros. Indeed setting $\mu=\nu$ in
Corollary~\ref{cr:phi-below-conc} we get the lower bound, while setting $\mu=\nu$ in
Corollary~\ref{cr:f-twosided} multiplying throughout by $f_{a,c}(\nu;x)^2$
and subtracting the same expression we get the upper bound.

\textbf{Remark.} Theorems~\ref{th:f-convexity} and \ref{th:g-convexity} are
independent in the sense that

\begin{theo}\label{th:g-concavity}
Suppose either \emph{(a)} $c+1\geq{a}\geq{c}>0$ and
$\{g_n\}_{n=0}^{\infty}$ is an arbitrary nonnegative sequence or
\emph{(b)} $a>c+1>1$ and $\{g_n\}_{n=0}^{\infty}$ is a nonnegative
log-concave sequence without internal zeros.
Then  $\mu\mapsto g_{a,c}(\mu;x)$ is discrete Wright
log-concave on $[0,\infty)$ for each fixed $x>0$.
Moreover, the family $\{g_{a,c}(\mu;x)\}_{\mu\geq0}$ is
discrete $q$-log-concave, i.e. the function $x\to\psi_{a,c}(\mu,\nu;x)$
has nonnegative power series coefficients for all $\nu\in\N$
and $\mu\geq\nu-1$ so that $x\to\psi_{a,c}(\mu,\nu;x)$ is absolutely
monotonic and multiplicatively convex on $(0,\infty)$.
\end{theo}

\textbf{Proof.} According to Lemmas~\ref{lm:disc-Wright} and \ref{lm:disc-qWright}
it is sufficient to consider the case $\nu=1$. For a fixed integer $m\geq{1}$ we have by the Cauchy
product and Gauss summation:
\begin{multline*}
\psi_m=\sum\limits_{k=0}^{m}\frac{g_{k}g_{m-k}}{k!(m-k)!}\left[\frac{\Gamma(a+1)(a+1)_k\Gamma(a+\mu)(a+\mu)_{m-k}}{\Gamma(c+1)(c+1)_k\Gamma(c+\mu)(c+\mu)_{m-k}}-
\frac{\Gamma(a)(a)_k\Gamma(a+\mu+1)(a+\mu+1)_{m-k}}{\Gamma(c)(c)_k\Gamma(c+\mu+1)(c+\mu+1)_{m-k}}\right]
\\
=\frac{\Gamma(a)\Gamma(a+\mu)}{\Gamma(c)\Gamma(c+\mu)}\sum\limits_{k=0}^{m}\frac{g_{k}g_{m-k}}{k!(m-k)!}\left[\frac{a(a+1)_k(a+\mu)_{m-k}}{c(c+1)_k(c+\mu)_{m-k}}-
\frac{(a)_k(a+\mu)(a+\mu+1)_{m-k}}{(c)_k(c+\mu)(c+\mu+1)_{m-k}}\right]
\\
=\frac{\Gamma(a)\Gamma(a+\mu)}{\Gamma(c)\Gamma(c+\mu)}
\sum\limits_{k=0}^{m}\frac{g_{k}g_{m-k}}{k!(m-k)!}\frac{(a)_{k}(a+\mu)_{m-k}}{(c)_{k}(c+\mu)_{m-k}}\left[\frac{(a+k)}{(c+k)}-
\frac{(a+\mu+m-k)}{(c+\mu+m-k)}\right]
\\
=\frac{\Gamma(a)\Gamma(a+\mu)}{\Gamma(c)\Gamma(c+\mu)}
\sum\limits_{k=0}^{m}\frac{g_{k}g_{m-k}}{k!(m-k)!}\frac{(a)_{k}(a+\mu)_{m-k}}{(c)_{k}(c+\mu)_{m-k}}\frac{(a-c)(m-2k+\mu)}{(c+k)(c+\mu+m-k)}
\\
=\frac{(a-c)\Gamma(a)\Gamma(a+\mu)}{\Gamma(c+1)\Gamma(c+\mu+1)m!}
\sum\limits_{k=0}^{m}\frac{g_{k}g_{m-k}(a)_{k}(a+\mu)_{m-k}}{(c+1)_{k}(c+1+\mu)_{m-k}}\binom{m}{k}(m-2k+\mu)
\\
=\frac{(a-c)\Gamma(a)\Gamma(a+\mu)}{\Gamma(c+1)\Gamma(c+\mu+1)m!}\sum\limits_{k=0}^{[m/2]}g_{k}g_{m-k}\binom{m}{k}M_k,
\end{multline*}
where
$$
M_k=\underbrace{\frac{(a)_{k}(a+\mu)_{m-k}}{(c+1)_{k}(c+1+\mu)_{m-k}}}_{=V_{k}}(m-2k+\mu)-\underbrace{\frac{(a)_{m-k}(a+\mu)_{k}}{(c+1)_{m-k}(c+1+\mu)_{k}}}_{=V_{m-k}}(m-2k-\mu).
$$
for $k<m/2$ and
$$
M_k=\frac{(a)_{k}(a+\mu)_{m-k}\mu}{(c+1)_{k}(c+1+\mu)_{m-k}}
$$
for $k=m/2$.  Lemma~\ref{lm:ab-sum} shows that
$$
\sum\limits_{0\leq{k}\leq{m/2}}\binom{m}{k}M_{k}>0
$$
for all $a>c>0$.
Moreover, the proof of the lemma for the case $c+1=b\geq{a}>0$
implies that $M_k>0$ for all $k=0,2,\ldots,[m/2]$.  This proves
the first part of the theorem. In order to prove the second part
pertaining to $a>c+1$ we will apply Lemma~\ref{lm:sum}.  Setting
$A_k=\binom{m}{k}M_k$ it is left to demonstrate that that sequence
$M_0,M_1,\ldots,M_{[m/2]}$ changes sign no more than once.
Indeed for $k=m-k$ it is clear that $M_k>0$. If $k<m-k$ then
$$
V_k<V_{m-k}~\Leftrightarrow~\frac{(a+\mu+k)\cdots(a+\mu+m-k-1)}{(c+1+\mu+k)\cdots(c+1+\mu+m-k-)}<\frac{(a+k)\cdots(a+m-k-1)}{(c+1+k)\cdots(c+1+m-k-1)}
$$
since $x\to(a+x)/(c+1+x)$ is decreasing.
Assume that $M_k<0$ for some $k$. We will demonstrate that
$M_{k-1}<0$.  We have
$$
M_{k-1}=V_kR(\mu)(m-2k+\mu+2)-V_{m-k}S(\mu)(m-2k-\mu+2),
$$
where
$$
R(\delta)=\frac{(a+\delta+m-k)(c+k)}{(c+1+\delta+m-k)(a+k-1)},~~
S(\delta)=\frac{(c+\delta+k)(a+m-k)}{(a-1+\delta+k)(c+1+m-k)}.
$$
It follows from $R(0)=S(0)$  and $V_{k}<V_{m-k}$ that
$$
V_kR(0)(m-2k+\mu+2)-V_{m-k}S(0)(m-2k-\mu+2)=R(0)(M_k+2(V_k-V_{m-k}))<0.
$$
Next, $R(\delta)$ is decreasing, while $S(\delta)$ is increasing
on $(0,\infty)$ because $a>c+1$, and $m-2k-\mu>0$ because  $M_k<0$
so that
$$
M_{k-1}=V_kR(\mu)(m-2k+\mu)-V_{m-k}S(\mu)(m-2k-\mu)+2(V_kR(\mu)-V_{m-k}S(\mu))
$$
$$
<R(0)(M_k+2(V_k-V_{m-k}))<0.~~~~\square
$$
Again we have the following corollaries.
\begin{corol}\label{cr:psi-complmon-conc}
Suppose $\nu\in\N$ and  $\mu\geq\nu-1$. Under hypotheses of
Theorem~\ref{th:g-concavity} the function $y\to\psi_{a,c}(\mu,\nu;1/y)$
is completely monotonic and log-convex on $[0,\infty)$
and there exists a nonnegative measure $\tau$ supported on $[0,\infty)$
such that
$$
\psi_{a,c}(\mu,\nu;x)=\int\limits_{[0,\infty)}e^{-t/x}d\tau(t).
$$
\end{corol}
\begin{corol}\label{cr:psi-below-conc}
Under hypotheses and notation of Theorem~\ref{th:g-concavity} and
for all $\nu\in\N$, $\mu\geq\nu-1$, $x\geq{0}$
$$
g_{a,c}(\mu;x)g_{a,c}(\nu;x)-g_{a,c}(0;x)g_{a,c}(\mu+\nu;x)
\!\geq\!
g_0^2\left\{\frac{\Gamma(a+\nu)\Gamma(a+\mu)}{\Gamma(c+\nu)\Gamma(c+\mu)}-
\frac{\Gamma(a)\Gamma(a+\mu+\nu)}{\Gamma(c)\Gamma(c+\mu+\nu)}\right\}
$$
with equality only at $x=0$ if $a-c,\mu,\nu\ne0$.
\end{corol}
\begin{corol}\label{cr:g-twosided}
Under hypotheses and notation of Theorem~\ref{th:g-concavity}
$$
\frac{(a+\mu)_{\nu}(c)_{\nu}}{(c+\mu)_{\nu}(a)_{\nu}}\leq\frac{g_{a,c}(0;x)g_{a,c}(\mu+\nu;x)}{g_{a,c}(\nu;x)g_{a,c}(\mu;x)}\leq1
$$
for all $\nu\in\N$, $\mu\geq{0}$ and $x\geq{0}$.
\end{corol}

Combining corollaries~\ref{cr:psi-below-conc} and \ref{cr:g-twosided} we obtain
the following two-sided bounds for the Tur\'{a}nian:
\begin{equation}\label{eq:g-Turanian}
g_0^2\frac{\Gamma(a)^2}{\Gamma(c)^2}\left[\frac{(a)_{\nu}^2}{(c)_{\nu}^2}-
\frac{(a)_{2\nu}}{(c)_{2\nu}}\right]\!\leq\!g_{a,c}(\nu;x)^2-g_{a,c}(0;x)g_{a,c}(2\nu;x)
\!\leq\!\frac{(c+\nu)_{\nu}(a)_{\nu}-(a+\nu)_{\nu}(c)_{\nu}}{(a)_{\nu}(c+\nu)_{\nu}}g_{a,c}(\nu;x)^2.
\end{equation}
The bounds are valid for $\nu\in\N$, $a\geq{c}>0$ and assuming that $\{g_n\}_{n\geq{0}}$ is
a nonnegative sequence which is also log-concave and without internal zeros if $a>c+1$.

There is virtually no doubt that the discrete Wright log-concavity
demonstrated in Theorem~\ref{th:g-concavity} results from our method
of proof so that the adjective ''discrete'' is redundant. In other
words, we propose the following conjecture.

\begin{conjecture} The family $\{g_{a,c}(\mu;x)\}_{\mu\geq{0}}$ is
Wright $q$-log-concave  for all $a\geq{c}>0$.
\end{conjecture}

The next theorem deals with the class of series  defined by
\begin{equation}\label{eq:h-def}
\mu\to
h_{a,c}(\mu;x)=\sum\limits_{n=0}^{\infty}h_n\frac{(a+\mu)_n}{\Gamma(c+\mu+n)}\frac{x^n}{n!}
\end{equation}
and their product differences
\begin{equation}\label{eq:lambda-def}
\lambda_{a,c}(\mu,\nu;x)=h_{a,c}(\mu;x)h_{a,c}(\nu;x)-h_{a,c}(0;x)h_{a,c}(\nu+\mu;x)
=\sum\limits_{m=0}^{\infty}\lambda_mx^m.
\end{equation}

Here we have discrete $q$-logarithmic concavity
 for all nonnegative values of $a$ and $c$.
\begin{theo}\label{th:h-concavity}
Suppose either \emph{(a)} $c+1\geq{a}\geq{c}>0$ and
$\{h_n\}_{n=0}^{\infty}$ is an arbitrary nonnegative sequence or
\emph{(b)}  $c>0$,  $a\in(0,c)\cup(c+1,\infty)$ and $\{h_n\}_{n=0}^{\infty}$
is a nonnegative log-concave sequence without internal zeros.
Then  $\mu\mapsto h_{a,c}(\mu;x)$ is discrete Wright
log-concave on $[0,\infty)$ for each fixed $x>0$.
Moreover, the family $\{h_{a,c}(\mu;x)\}_{\mu\geq0}$ is
discrete $q$-log-concave, i.e. $x\to\lambda_{a,c}(\mu,\nu;x)$
has nonnegative power series coefficients for $\nu\in\N$,
$\mu\geq{\nu-1}$ so that $x\to\lambda_{a,c}(\mu,\nu;x)$ is
absolutely monotonic and multiplicatively convex on $(0,\infty)$.
\end{theo}
\textbf{Proof.} Suppose $a\geq{c}>0$. We have
$$
h_{a,c}(\mu;x)=\frac{1}{\Gamma(a+\mu)}g_{a,c}(\mu;x),
$$
where $g_{a,c}(\mu;x)$ is defined in (\ref{eq:g-def}) and $g_n=h_n$.
The claims of the theorem for $a\geq{c}>0$ then follow from Theorem~\ref{th:g-concavity}
and Lemma~\ref{lm:product}.

If $c>a>0$ write
$$
h_{a,c}(\mu;x)=\frac{1}{\Gamma(c+\mu)}f_{a,c}(\mu;x),
$$
where $f_{a,c}(\mu;x)$ is defined in (\ref{eq:f-def}) and $f_n=h_n$.
The claims of the theorem for $c>a>0$ then follow from Theorem~\ref{th:f-concavity}
and Lemma~\ref{lm:product}.~$\square$

\begin{corol}\label{cr:lambda-complmon-conc}
Under hypotheses of Theorem~\ref{th:h-concavity}
and assuming the series in $(\ref{eq:h-def})$ converges for all $x\geq{0}$
the function $y\to\lambda_{a,c}(\mu,\nu;1/y)$ is completely
monotonic and log-convex on $[0,\infty)$ for all $\nu\in\N$ and
$\mu\geq\nu-1$ so that there exists a nonnegative measure $\tau$
supported on $[0,\infty)$ such that
$$
\lambda_{a,c}(\mu,\nu;x)=\int\limits_{[0,\infty)}e^{-t/x}d\tau(t).
$$
\end{corol}
\begin{corol}\label{cr:lambda-below-conc}
Under hypotheses of Theorem~\ref{th:h-concavity} and
for all $\nu\in\N$, $\mu\geq{\nu-1}$ and $x\geq{0}$
$$
h_{a,c}(\mu;x)h_{a,c}(\nu;x)-h_{a,c}(0;x)h_{a,c}(\mu+\nu;x)\geq
\frac{h_0^2[(c+\mu)_{\nu}-(c)_{\nu}]}{\Gamma(c+\nu)\Gamma(c+\mu+\nu)}
$$
with equality only at $x=0$ if $\mu,\nu\ne0$.
\end{corol}

Finally, we consider the the class of series defined by
\begin{equation}\label{eq:q-def}
\mu\to
q_{a,c}(\mu;x)=\sum\limits_{n=0}^{\infty}q_n\frac{\Gamma(a+\mu+n)}{(c+\mu)_n}\frac{x^n}{n!}
\end{equation}
and their product differences
\begin{equation}\label{eq:rho-def}
\rho_{a,c}(\mu,\nu;x)=q_{a,c}(\mu;x)q_{a,c}(\nu;x)-q_{a,c}(0;x)q_{a,c}(\nu+\mu;x)
=\sum\limits_{m=0}^{\infty}\rho_mx^m.
\end{equation}

Here we have $q$-logarithmic convexity for all nonnegative values of $a$ and $c$.
\begin{theo}\label{th:q-convexity}
Suppose $a,c>0$, $\{q_n\}_{n=0}^{\infty}$ is any nonnegative
sequence and the functions $q_{a,c}(\mu;x)$ and  $\rho_{a,c}(\mu,\nu;x)$
are defined by $(\ref{eq:q-def})$ and $(\ref{eq:rho-def})$, respectively.
Then the function $\mu\mapsto q_{a,c}(\mu;x)$ is strictly
log-convex on $[0,\infty)$ for each fixed $x>0$. Moreover,
the family $\{q_{a,c}(\mu;x)\}_{\mu\geq0}$ is
Wright $q$-log-concave, i.e. the
function $x\to\rho_{a,c}(\mu,\nu;x)$ has non-positive  power
series coefficients for all $\mu,\nu\geq{0}$ so that the function
$x\to-\rho_{a,c}(\mu,\nu;x)$ is absolutely monotonic and
multiplicatively convex on
$(0,\infty)$.
\end{theo}
\textbf{Proof.} Suppose $a\geq{c}>0$. We have
$$
q_{a,c}(\mu;x)=\Gamma(a+\mu)f_{a,c}(\mu;x),
$$
where $f_{a,c}(\mu;x)$ is defined in (\ref{eq:f-def}) and $f_n=q_n$.
The claims of the theorem for $a\geq{c}>0$ then follow from Theorem~\ref{th:f-convexity}
and Lemma~\ref{lm:product}.

If $c>a>0$ write
$$
q_{a,c}(\mu;x)=\Gamma(c+\mu)g_{a,c}(\mu;x),
$$
where $g_{a,c}(\mu;x)$ is defined in (\ref{eq:g-def}) and $g_n=q_n$.
The claims of the theorem for $c>a>0$ then follow from Theorem~\ref{th:g-convexity}
and Lemma~\ref{lm:product}.~$\square$

\begin{corol}\label{cr:rho-below-conv}
Under hypotheses of Theorem~\ref{th:q-convexity} and
for all $x\geq{0}$
$$
q_{a,c}(0;x)q_{a,c}(\mu+\nu;x)-q_{a,c}(\mu;x)q_{a,c}(\nu;x)
\!\geq\!
q_0^2\left\{\Gamma(a)\Gamma(a+\mu+\nu)-\Gamma(a+\mu)\Gamma(a+\nu)\right\}
$$
with equality only at $x=0$  if $\mu,\nu\ne0$.
\end{corol}
\begin{corol}\label{cr:rho-complmon-conv}
Under hypotheses of Theorem~\ref{th:q-convexity} and assuming the series
in $(\ref{eq:q-def})$ converges for all $x\geq{0}$ the function
$y\to-\rho_{a,c}(\mu,\nu;1/y)$
is completely monotonic and log-convex on $[0,\infty)$ for all $\mu,\nu>0$
and there exists a nonnegative measure $\tau$ supported on $[0,\infty)$
such that
$$
\rho_{a,c}(\mu,\nu;x)=-\int\limits_{[0,\infty)}e^{-t/x}d\tau(t).
$$
\end{corol}

\begin{corol}\label{cr:h-twosided}
Under hypotheses and notation of Theorem~\ref{th:h-concavity}
$$
\frac{(a)_{\nu}(c)_{\nu}}{(a+\mu)_{\nu}(c+\mu)_{\nu}}\leq\frac{h_{a,c}(0;x)h_{a,c}(\mu+\nu;x)}{h_{a,c}(\nu;x)h_{a,c}(\mu;x)}\leq1
$$
for $\nu\in\N$, $x,\mu\geq{0}$.
\end{corol}
Combining corollaries~\ref{cr:lambda-below-conc} and \ref{cr:h-twosided} we obtain
the following two-sided bounds for the Tur\'{a}nian:
\begin{equation}\label{eq:h-Turanian}
\frac{h_0^2[(c+\nu)_{\nu}-(c)_{\nu}]}{\Gamma(c+\nu)\Gamma(c+2\nu)}\leq{h_{a,c}(\nu;x)^2-h_{a,c}(0;x)h_{a,c}(2\nu;x)}\leq\left(1-\frac{(a)_{\nu}(c)_{\nu}}{(a+\nu)_{\nu}(c+\nu)_{\nu}}\right)h_{a,c}(\nu;x)^2.
\end{equation}
The bounds are valid for $\nu\in\N$, $a,c>0$ and assuming that $\{h_n\}_{n\geq{0}}$ is
a nonnegative sequence which is also log-concave and without internal zeros if
$a\in(0,c)\cup(c+1,\infty)$.

\bigskip

\paragraph{4. Applications and relation to other work.}
In this section we will demonstrate how Theorems~\ref{th:f-concavity}
to \ref{th:q-convexity} and their corollaries lead to various
new inequalities for the Kummer, Gauss and generalized hypergeometric
functions and their ratios and logarithmic derivatives.

\medskip

\textbf{Example~1.} The first very natural candidate to apply the theory
presented in the previous section is the Kummer function ${_1F_1}(a;c;x)$.
Indeed, setting
$$
f_{a,c}(\mu;x)={_1F_1}(a+\mu;c+\mu;x),~~~~g_{a,c}(\mu;x)=
\frac{\Gamma(a+\mu)}{\Gamma(c+\mu)}{_1F_1}(a+\mu;c+\mu;x),
$$
$$
h_{a,c}(\mu;x)=\frac{1}{\Gamma(c+\mu)}{_1F_1}(a+\mu;c+\mu;x)~~\text{and}~~q_{a,c}(\mu;x)=
\Gamma(a+\mu){_1F_1}(a+\mu;c+\mu;x)
$$
we obtain examples of functions defined
by (\ref{eq:f-def}), (\ref{eq:g-def}),
(\ref{eq:h-def}) and (\ref{eq:q-def}),
respectively, and satisfying the corresponding
theorems and corollaries.  These facts extend and refine
some previous results due to Baricz \cite[Theorem~2]{BariczExpo} and
the second author \cite{KarpPOMI}.  In particular, we obtain the following
bounds for the Tur\'{a}nian:
\begin{equation}\label{eq:Turan1}
\frac{2x(c-a)}{(c)_3}\leq {_1F_1}(a+1;c+1;x)^2-{_1F_1}(a;c;x){_1F_1}(a+2;c+2;x)
\leq\frac{c-a}{c(a+1)}{_1F_1}(a+1;c+1;x)^2
\end{equation}
if $c\geq{a}>0$ and
\begin{equation}\label{eq:Turan2}
\frac{a-c}{c(c+1)}
\!\leq\!\frac{a}{c}{_1F_1}(a+1;c+1;x)^2
-\frac{a+1}{c+1}{_1F_1}(a;c;x){_1F_1}(a+2;c+2;x)
\!\leq\!\frac{a-c}{c(c+1)}{_1F_1}(a+1;c+1;x)^2
\end{equation}
if $a\geq{c}>0$.  Simple rearrangements of the
righthand side of (\ref{eq:Turan1}) and the lefthand side
of (\ref{eq:Turan2}) give
\begin{equation}\label{eq:Turan3}
\frac{a}{c}{_1F_1}(a+1;c+1;x)^2
-\frac{a+1}{c+1}{_1F_1}(a;c;x){_1F_1}(a+2;c+2;x)~\left\{
\!\!\!\begin{array}{l}
\leq0,~~c\geq{a}>0,\\[5pt]
\geq0,~~a\geq{c}>0.
\end{array}\right.
\end{equation}
Substituting the contiguous relation
\begin{equation}\label{eq:contiguous}
{_1F_1}(a+2;c+2;x)=\frac{(c+1)(x-c)}{(a+1)x}{_1F_1}(a+1;c+1;x)
+\frac{c(c+1)}{(a+1)x}{_1F_1}(a;c;x)
\end{equation}
into (\ref{eq:Turan3}) we get after some algebra:
$$
xy^2+(c-x)y-a~\left\{
\!\!\!\begin{array}{l}
\leq0,~~c\geq{a}>0,\\[5pt]
\geq0,~~a\geq{c}>0.
\end{array}\right.
$$
where
$$
y=\frac{{_1F_1}\!'(a;c;x)}{{_1F_1}(a;c;x)}=\frac{a{_1F_1}(a+1;c+1;x)}{c{_1F_1}(a;c;x)}.
$$
In a similar fashion writing (\ref{eq:Turan3})
with $a\to{a+1}$, $c\to{c+1}$ we obtain:
$$
\frac{a+1}{c+1}{_1F_1}(a+2;c+2;x)^2-\frac{a+2}{c+2}{_1F_1}(a+1;c+1;x){_1F_1}(a+3;c+3;x)~\left\{
\!\!\!\begin{array}{l}
\leq0,~~c\geq{a}>0,\\[5pt]
\geq0,~~a\geq{c}>0.
\end{array}\right.
$$
Using contiguous relation (\ref{eq:contiguous}) twice this leads to
$$
(ax+c)y^2-a(x-c+1)y-a^2~\left\{
\!\!\!\begin{array}{l}
\geq0,~~c\geq{a}>0,\\[5pt]
\leq0,~~a\geq{c}>0.
\end{array}\right.
$$
Solving the two quadratics we obtain
$$
\frac{x-c+1+\sqrt{(x-c+1)^2+4ax+4c}}{2x+2c/a}\leq
\frac{{_1F_1}\!'(a;c;x)}{{_1F_1}(a;c;x)}\leq\frac{x-c+\sqrt{(x-c)^2+4ax}}{2x}
$$
if $c\geq{a}>0$ and
$$
\frac{x-c+\sqrt{(x-c)^2+4ax}}{2x}\leq
\frac{{_1F_1}\!'(a;c;x)}{{_1F_1}(a;c;x)}\leq\frac{x-c+1+\sqrt{(x-c+1)^2+4ax+4c}}{2x+2c/a}
$$
if $a\geq{c}>0$. Note that for $a=c$ both bounds reduce to $1$.
It is also easy to check that both bounds give correct value $1$ at
$x=\infty$ and correct value $a/c$ at $x=0$.  Moreover,
the upper bound in the first inequality has a correct term of
order $O(1/x)$ around infinity, while the lower bound has a correct term
of order $O(x)$ around zero.  Note that similar
but different bounds have been obtained in our recent paper \cite{KK}.

Integrating these bound from $0$ to $x$ we obtain
$$
B_1(x)\leq{_1F_1}(a;c;x)\leq{B_2(x)}~\text{if}~c\geq{a}>0;~~B_2(x)\leq{_1F_1}(a;c;x)\leq{B_1(x)}~\text{if}~a\geq{c}>0,
$$
where (we set $b=(a+1)(a-c)$ for brevity)
$$
B_1(x)={\small
\frac{(2+2a)^{-b/a}c^{(a^2-b)/a}\bigl(1+2a-c+x+\sqrt{(x-c+1)^2+4ax+4c}\bigr)^{(a^2+b)/2a}}
{\bigl(c^2(a+1)+a-c+(a^2+b)x+(a^2-b)\sqrt{(x-c+1)^2+4ax+4c}\bigr)^{(a^2-b)/2a}}\times}
$$
$$
\times \exp\left\{\dfrac{x-c-1+\sqrt{(x-c+1)^2+4ax+4c}}{2}\right\}
$$
$$
B_2(x)=\frac{(4ac)^{c/2}}{(2a)^a}\frac{(2a+\sqrt{(x-c)^2+4ax}+x-c)^{a-c/2}}{(2ax/c+\sqrt{(x-c)^2+4ax}-(x-c))^{c/2}}\exp\left\{\frac{\sqrt{(x-c)^2+4ax}+x-c}{2}\right\}
$$
All the above bounds can be easily extended to
$x<0$ using the Kummer identity ${_1F_1}(a;c;x)=e^x{_1F_1}(c-a;c;-x)$.

\textbf{Example~2.}  The ratio
$$
r(x):=\frac{{_2F_1}(a+1,b;c+1;x)}{{_2F_1}(a,b;c;x)}
$$
was first developed into continued fraction by Euler. Later, Gauss found
a different continued fraction which became more popular than the original fraction
of Euler, see \cite[paragraph 2.5]{AAR} for details and references.  Here we will derive bounds for this ratio under some restrictions on
parameters which are related in some way (see below) to
the continued fraction of Euler.  According to \cite[formula (2.5.3)]{AAR}
\begin{equation}\label{eq:2F1contig}
\frac{a+1}{c+1}{_2F_1}(a+2,b;c+2;x)=\frac{c+(a-b+1)x}{(c-b+1)x}{_2F_1}(a+1,b;c+1;x)
-\frac{c}{(c-b+1)x}{_2F_1}(a,b;c;x).
\end{equation}
Further, setting  $g_n=(b)_n$ in (\ref{eq:g-def}) we get
$$
g_{a,c}(\mu;x)=\frac{\Gamma(a+\mu)}{\Gamma(c+\mu)}{_2F_1}(a+\mu,b;c+\mu;x).
$$
Then it follows from Corollary~\ref{cr:psi-below-conv} (with $\mu=\nu=1$ and using $g_0>0$)
that
\begin{equation}\label{eq:2F1Turan}
\frac{a}{c}({_2F_1}(a+1,b;c+1;x))^2\leq\frac{a+1}{c+1}{_2F_1}(a,b;c;x){_2F_1}(a+2,b;c+2;x),~~~0\leq{x}<1,
\end{equation}
if $c\geq{a}>0$.  Substituting (\ref{eq:2F1contig}) here we obtain
$$
\frac{a}{c}({_2F_1}(a+1,b;c+1;x))^2\leq\frac{c+(a-b+1)x}{(c-b+1)x}{_2F_1}(a,b;c;x){_2F_1}(a+1,b;c+1;x)
-\frac{c({_2F_1}(a,b;c;x))^2}{(c-b+1)x}
$$
or, after division by $({_2F_1}(a,b;c;x))^2$,
$$
\frac{a}{c}r(x)^2-\frac{c+(a-b+1)x}{(c-b+1)x}r(x)+\frac{c}{(c-b+1)x}\leq0þ
$$
Solving this quadratic inequality for $c-b+1<0$ and $c-b+1>0$ we arrive at
$$
r(x)\leq\frac{c+(a-b+1)x-\sqrt{(c+(a-b+1)x)^2-4a(c-b+1)x}}{2(a/c)(c-b+1)x},~\text{if}~c+1<b,
$$
$$
r(x)\geq\frac{c+(a-b+1)x-\sqrt{(c+(a-b+1)x)^2-4a(c-b+1)x}}{2(a/c)(c-b+1)x},~\text{if}~c+1>b.
$$
Note that for $c=b+1$ both inequalities turn into correct equality $r(x)=c/(c-(c-a)x)$.
It is also easy to verify that $r(0)=1$ coincides with the value of the bound at $x=0$.
Using rather standard techniques the expression on the right of the two formulas
above can be developed into continued fraction:
$$
\frac{c+(a-b+1)x-\sqrt{(c+(a-b+1)x)^2-4a(c-b+1)x}}{2(a/c)(c-b+1)x}=\frac{c}{a(b-c-1)x}
\K\limits_{n=0}^{\infty}\frac{a(b-c-1)x}{c+(a-b+1)x}
$$
which is interesting to compare with the continued fraction of Euler:
$$
r(x)=\frac{c}{a(b-c)x}\K\limits_{n=0}^{\infty}\frac{(a+n)(b-c-n)x}{c+n+(a-b+n+1)x}.
$$
Here, we employed the usual notation
$$
\K\limits_{n=0}^{\infty}\frac{a_n}{b_n}=\cfrac{a_0}{b_0+\cfrac{a_1}{b_1+\cdots}}.
$$
If the last fraction for $r(x)$ is made 1-periodic starting from $n=0$,
$$
\frac{c}{a(b-c)x}\K\limits_{n=0}^{\infty}\frac{a(b-c)x}{c+(a-b+1)x}
=\frac{c+(a-b+1)x-\sqrt{(c+(a-b+1)x)^2-4a(c-b)x}}{2(a/c)(c-b)x},
$$
we get an approximation which, by numerical tests, underestimates $r(x)$ and is less precise
then our bound above.  We can obtain a sequence of improving approximations to  $r(x)$
by continued fractions which are 1-periodic starting from $n=N$, $N=1,2,\ldots$.
Each approximation in this sequence is a rational function of $x$ and
square root of some quadratic of $x$ and is easily computable.

We note the the above bounds can be extended to negative values of the
argument by an application of Pfaff's transformation \cite[formula (2.2.6)]{AAR}
$$
{_2F_1}(a,b;c;x)=(1-x)^{-a}{_2F_1}(a,c-b;c;x/(x-1)).
$$

\textbf{Example~3.}
The application of Theorems~\ref{th:f-concavity} to \ref{th:q-convexity}
to generalized hypergeometric function is largely based on the following
lemma.
\begin{lemma}\label{lm:HVVKS}
Denote by $e_k(x_1,\ldots,x_q)$  the $k$-th elementary symmetric
polynomial,
$$
e_0(x_1,\ldots,x_q)=1,~~~e_k(x_1,\ldots,x_q)=\!\!\!\!\!\!\!\!\sum\limits_{1\leq{j_1}<{j_2}\cdots<{j_k}\leq{q}}
\!\!\!\!\!\!\!\!x_{j_1}x_{j_2}\cdots{x_{j_k}},~~k\geq{1}.
$$
Suppose $q\geq{1}$ and $0\leq{r}\leq{q}$ are integers, $a_i>0$,
$i=1,\ldots,q-r$, $b_i>0$, $i=1,\ldots,q$, and
\begin{equation}\label{eq:symmetric-chain1}
\frac{e_q(b_1,\ldots,b_q)}{e_{q-r}(a_1,\ldots,a_{q-r})}\leq
\frac{e_{q-1}(b_1,\ldots,b_q)}{e_{q-r-1}(a_1,\ldots,a_{q-r})}\leq\cdots
\leq\frac{e_{r+1}(b_1,\ldots,b_q)}{e_{1}(a_1,\ldots,a_{q-r})}\leq
e_{r}(b_1,\ldots,b_q).
\end{equation}
Then the sequence of hypergeometric terms \emph{(}if $r=q$ the
numerator is $1$\emph{)},
$$
f_n=\frac{(a_1)_n\cdots(a_{q-r})_n}{(b_1)_n\cdots(b_q)_n},
$$
is log-concave, i.e. $f_{n-1}f_{n+1}\leq{f_n^2}$, $n=1,2,\ldots$
It is strictly log-concave unless $r=0$ and $a_i=b_i$,
$i=1,\ldots,q$.
\end{lemma}
The proof of this lemma for $r=0$ follows from
\cite[Theorem~4.4]{HVV}. For general $r$ see \cite[Lemma~2]{KS}
and the last paragraph of that paper.

We note that it has been demonstrated in \cite[Lemma~2]{KarpPOMINew} that
(\ref{eq:symmetric-chain1}) is true for $r=0$ if majorization conditions
$$
\sum\limits_{i=1}^{k}b_i\leq\sum\limits_{i=1}^{k}a_i~~\text{for}~~k=1,2,\ldots,q,
$$
hold, where
$$
0<a_1\leq{a_2}\leq\cdots\leq{a_q},~~0<b_1\leq{b_2}\leq\cdots\leq{b_q}.
$$

Consider the following functions  ($p,q\geq{1}$):
$$
f_{a,c}(\mu;x)={_pF_q}(a+\mu,a_2,\ldots,a_p;c+\mu,c_2,\ldots,c_q;x),
$$
$$
g_{a,c}(\mu;x)=
\frac{\Gamma(a+\mu)}{\Gamma(c+\mu)}{_pF_q}(a+\mu,a_2,\ldots,a_p;c+\mu,c_2,\ldots,c_q;x),
$$
$$
h_{a,c}(\mu;x)=\frac{1}{\Gamma(c+\mu)}{_pF_q}(a+\mu,a_2,\ldots,a_p;c+\mu,c_2,\ldots,c_q;x)
$$
and
$$
q_{a,c}(\mu;x)=\Gamma(a+\mu){_pF_q}(a+\mu,a_2,\ldots,a_p;c+\mu,c_2,\ldots,c_q;x).
$$
Assuming all parameters are positive these functions satisfy Theorems~\ref{th:f-convexity}, \ref{th:g-convexity},
\ref{th:g-concavity}(a), \ref{th:h-concavity}(a), \ref{th:q-convexity} and their corollaries
without any further restrictions. If, in addition, $p\leq{q}$ and the
vectors $(a_2,\ldots,a_p)$ and $(b_2,\ldots,b_q)$ satisfy Lemma~\ref{lm:HVVKS}
with $r=q-p$ then $f_{a,c}(\mu;x)$, $g_{a,c}(\mu;x)$ and $h_{a,c}(\mu;x)$
satisfy Theorems~\ref{th:f-concavity}, \ref{th:g-concavity}(b) and \ref{th:h-concavity}(b),
respectively.  These facts imply a number of presumably new inequalities for the generalized
hypergeometric function.  In particular, if $\nu\in\N$, $x\geq{0}$ and under conditions (\ref{eq:symmetric-chain1})
the function $f_{a,c}(\mu;x)$ defined above satisfies (\ref{eq:f-Turanian}) for $c\geq{a}>0$,
$g_{a,c}(\mu;x)$ satisfies (\ref{eq:g-Turanian}) for $a\geq{c}>0$ and
$h_{a,c}(\mu;x)$ satisfies (\ref{eq:h-Turanian}) for all $a,c>0$.

\paragraph{5. Acknowledgements.} We thank Sergei Sitnik for
reading the manuscript and several useful remarks.
The research was supported by the Ministry of Education and
Science of Russia, projects 14.A18.21.0353 and 1.1773.2011, and the
Russian Foundation for Basic Research, project 11-01-00038-a.

\end{document}